\documentclass[UTF-8,reqno]{amsart}
\usepackage{enumerate}
\usepackage{mhequ}
\usepackage{amssymb,url,color, booktabs,nccmath}
\usepackage[left=3.2cm,right=3.2cm,top=4cm,bottom=4cm]{geometry}
\usepackage{mathrsfs}
\usepackage{enumitem,dsfont}

\usepackage{color}
\usepackage[colorlinks=true]{hyperref}
\hypersetup{
    linkcolor=blue,          
    citecolor=red,        
    filecolor=blue,      
    urlcolor=cyan
}

\definecolor{darkergreen}{rgb}{0.0, 0.5, 0.0}

\numberwithin{equation}{section}
\def\theequation{\arabic{section}.\arabic{equation}}
\newcommand{\be}{\begin{eqnarray}}
\newcommand{\ee}{\end{eqnarray}}
\newcommand{\ce}{\begin{eqnarray*}}
\newcommand{\de}{\end{eqnarray*}}
\newtheorem{theorem}{Theorem}[section]
\newtheorem{lemma}[theorem]{Lemma}
\newtheorem{proposition}[theorem]{Proposition}
\newtheorem{Examples}[theorem]{Example}
\newtheorem{corollary}[theorem]{Corollary}

\newtheorem{definition}[theorem]{Definition}
\theoremstyle{definition}
\newtheorem{remark}[theorem]{Remark}

\newcommand{\rmb}[1]{\textcolor{black}{#1}}

\newcommand{\tmop}[1]{\ensuremath{\operatorname{#1}}}

\def\${|\!|\!|}

\DeclareMathOperator{\supp}{supp}

\def\eps{\varepsilon}

\def\p{\partial}

\def\<{{\langle}}
\def\>{{\rangle}}
\def\({{\Big(}}
\def\){{\Big)}}

\def\bx{{\mathbf{x}}}
\def\tr{\mathrm {tr}}

\def\dif{{\mathord{{\rm d}}}}

\def\min{{\mathord{{\rm min}}}}

\def\no{\nonumber}
\def\={&\!\!=\!\!&}

\def\cB{{\mathcal B}}

\def\cF{{\mathcal F}}

\def\cL{{\mathcal L}}

\def\cP{{\mathcal P}}

\def\cR{{\mathcal R}}
\def\cS{{\mathcal S}}
\def\cT{{\mathcal T}}

\def\mN{{\mathbb N}}

\def\mP{{\mathbb P}}

\def\mR{{\mathbb R}}

\def\mT{{\mathbb T}}

\def\bP{{\mathbf P}}

\def\1{{\mathbf{1}}}

\def\E{\mathbf E}

\def\geq{\geqslant}
\def\leq{\leqslant}

\def\div{\mathord{{\rm div}}}

\def\eps{\varepsilon}

\def\p{\partial}

\def\<{{\langle}}
\def\>{{\rangle}}
\def\({{\Big(}}
\def\){{\Big)}}

\def\bx{{\mathbf{x}}}
\def\tr{\mathrm {Tr}}

\def\dif{{\mathord{{\rm d}}}}

\def\min{{\mathord{{\rm min}}}}

\def\no{\nonumber}
\def\={&\!\!=\!\!&}
\def\bt{\begin{theorem}}
\def\et{\end{theorem}}
\def\bl{\begin{lemma}}
\def\el{\end{lemma}}
\def\br{\begin{remark}}
\def\er{\end{remark}}
\def\bx{\begin{Examples}}
\def\ex{\end{Examples}}
\def\bd{\begin{definition}}
\def\ed{\end{definition}}
\def\bp{\begin{proposition}}
\def\ep{\end{proposition}}
\def\bc{\begin{corollary}}
\def\ec{\end{corollary}}

\def\geq{\geqslant}
\def\leq{\leqslant}

\def\div{\mathord{{\rm div}}}

\def\Id{\textrm{Id}}

\def\bP{{\mathbf P}}

 \def\R{\mathbb R}
 \def\R{\mathbb R}    
\def\N{\mathbb N}  
   
\def\<{\langle} \def\>{\rangle}

\allowdisplaybreaks

\begin{document}

\title[Non-unique ergodicity for 3D Navier--Stokes and Euler equations]{Non-unique ergodicity for deterministic and stochastic 3D Navier--Stokes and Euler equations}

\author{Martina Hofmanov\'a}
\address[M. Hofmanov\'a]{Fakult\"at f\"ur Mathematik, Universit\"at Bielefeld, D-33501 Bielefeld, Germany}
\email{hofmanova@math.uni-bielefeld.de}

\author{Rongchan Zhu}
\address[R. Zhu]{Department of Mathematics, Beijing Institute of Technology, Beijing 100081, China}
\email{zhurongchan@126.com}

\author{Xiangchan Zhu}
\address[X. Zhu]{ Academy of Mathematics and Systems Science,
Chinese Academy of Sciences, Beijing 100190, China}
\email{zhuxiangchan@126.com}
\thanks{
M.H. is grateful for funding from the European Research Council (ERC) under the European Union's Horizon 2020 research and innovation programme (grant agreement No. 949981) and  	the financial supports  by the Deutsche Forschungsgemeinschaft (DFG, German Research Foundation) – Project-ID 317210226--SFB 1283.
R.Z. and X.Z. are grateful to
the financial supports   by National Key R\&D Program of China (No. 2022YFA1006300).
R.Z. gratefully acknowledges financial support from the NSFC (No.
11922103, 12271030).
  X.Z. is grateful to the financial supports  in part by National Key R\&D Program of China (No. 2020YFA0712700) and the NSFC (No.   12090014, 12288201) and
  the support by key Lab of Random Complex Structures and Data Science,
 Youth Innovation Promotion Association (2020003), Chinese Academy of Science.
}

\begin{abstract}
We establish the existence of infinitely many stationary solutions, as well as  ergodic stationary solutions, to the three dimensional Navier--Stokes and Euler equations in both  deterministic and stochastic settings, driven by additive noise. These solutions belong to the regularity class $C(\mR;H^{\vartheta})\cap C^{\vartheta}(\mR;L^{2})$ for some $\vartheta>0$ and satisfy the equations in an analytically weak sense. The solutions to the Euler equations  are obtained as vanishing viscosity limits of stationary solutions to the Navier--Stokes equations. Furthermore, regardless of their construction, every stationary solution to the Euler equations within this regularity class, which satisfies  a suitable moment bound, is  a limit in law of stationary analytically weak solutions to Navier--Stokes equations with vanishing viscosities.
Our results are based on a novel stochastic version of the convex integration method, which provides uniform moment bounds locally in the aforementioned function spaces.
\end{abstract}


\date{\today}

\maketitle

\tableofcontents

\section{Introduction}

\subsection{Motivation}

Hydrodynamic turbulence is omnipresent in engineering applications and natural phenomena. Yet, developing a rigorous mathematical understanding of turbulence remains one of the major challenges in contemporary fluid dynamics research. To date,  results providing reliable predictions are very limited. On the physical side, the understanding has been driven by well-accepted theoretical hypotheses, such as those of the celebrated Kolmogorov theory  \cite{Kol41a, Kol41b, Kol41c}, see also \cite{Fri95}. These hypotheses have been largely confirmed by experiments. However, their rigorous verification from basic physical principles, particularly from the  incompressible Navier--Stokes equations, remains an outstanding open problem.
We encourage the reader to consult e.g. \cite{Fri95,MY13} for
thorough expositions of this topic.

The Navier--Stokes equations describe the time evolution of the velocity $u:[0,\infty)\times D\to \R^{3}$ of a viscous fluid confined in a domain $D\subset \R^{3}$. They are formulated  as follows:
\begin{align}\label{eq:detns}
\begin{aligned}
\partial_{t}u+\div(u\otimes u)+\nabla P&=\nu\Delta u +f,\\
\div u&=0,
\end{aligned}
\end{align}
where $P:[0,\infty)\times D\to \R$ denotes the associated pressure, $\nu>0$ the kinematic viscosity of the fluid and $f:[0,\infty)\times D\to \R^{3}$ is a given external force. These equations are further supplemented by initial and boundary conditions. Particularly relevant for the study of turbulence is the regime of high Reynolds number which corresponds to the vanishing viscosity limit $\nu\to0$. Formally, the Navier--Stokes equations then converge to the Euler equations:
\begin{align}\label{eq:deteu}
\begin{aligned}
\partial_{t}u+\div(u\otimes u)+\nabla P&=f,\\
\div u&=0,
\end{aligned}
\end{align}
which represent an idealized model for the highly turbulent limit regime. Assumptions allowing for a rigorous passage to the limit $\nu\to0$ are predicted by  physical theories of turbulence. However, it was shown in \cite{CG12} and \cite{CV18} that even under  weaker assumptions, convergence of the Navier--Stokes to Euler equations can be achieved, thus establishing the vanishing viscosity limit.

Experimental evidence has shown that  exact realizations of turbulent  trajectories are not suited for predictions due to their high sensitivity to input data such as initial and boundary conditions. On the other hand, and rather surprisingly, statistical properties are universal and well-reproducible.
Thus, a  probabilistic description appears indispensable. This necessity is further underlined by the lack of uniqueness, recently established for both the Euler and Navier--Stokes equations in the deterministic and stochastic settings, see e.g. \cite{BCV18, BMS20, BV19a,   CL20, DK20, DelSze2, DelSze3,  DelSze13,  L19, HZZ19, HZZ20, HZZ21, HZZ21markov}.

Furthermore, one of the fundamental assumptions in turbulence theory is the so-called ergodic hypothesis, which is taken for granted by physicists and engineers. This hypothesis asserts that time averages along trajectories coincide with ensemble averages taken with respect to some probability measure. This measure  is then invariant, i.e., preserved by the flow. Consequently, statistically stationary solutions (i.e., solutions whose probability law does not change over time)   play a distinguished role in the modeling of turbulence. It is  of essential interest to characterize these solutions and
 their attraction properties.

Moreover, it may be  possible to profit from the presence of stochastic perturbations of the equations, as some properties of the Navier--Stokes system have indeed been shown to improve under  stochastic noise. Specifically, the force $f$ is considered to be  a Gaussian noise that is white in time and colored in space, modeling large scale stirring that drives  turbulent fluids. In the deterministic setting, a selection of solutions depending continuously on the initial condition has not been obtained. However,  the probabilistic counterpart, i.e., the Feller property and even the strong Feller property which corresponds to a smoothing with respect to the initial condition, was established for a sufficiently non-degenerate noise in \cite{DD03} and \cite{FR08}. This led, in particular, to uniqueness of the invariant measure associated to the constructed selection of a Markov semigroup in these works.
Additionally, it was demonstrated in \cite{CH21} that   non-unique ergodic measures can be constructed in the Lorenz system by adding 
	noise to the last component. In \cite[Remark 1.3]{CH21}, the authors also pose the question of    whether a
	bifurcation of invariant measures appears at high Reynolds number for the Navier--Stokes system.

\medskip

It is therefore desired to investigate the validity of the following claims.
\begin{enumerate}
\item[(i)] The existence and (non)uniqueness of ergodic stationary solutions $u_{\nu}$ to the Navier--Stokes equations \eqref{eq:detns}.
\item[(ii)] The relative compactness of the family of stationary solutions $u_{\nu}$, $\nu>0$, and their convergence towards a statistically stationary solution to the Euler equations \eqref{eq:deteu}.
\item[(iii)] The existence and (non)uniqueness of ergodic stationary solutions to the Euler equations \eqref{eq:deteu}.
\end{enumerate}
To date,  these questions have only been addressed in several simplified settings, such as  certain shell models of turbulence \cite{FGHV14} and passive scalar models of turbulence \cite{BBPS19}. Additionally,  stationary solutions to 2D deterministic Euler equations, as limits of 2D stochastic Navier--Stokes equations, were obtained in \cite{K04}, and for the stochastic 1D Burgers equation \cite{BK21}. However, for the actual models of interest, i.e., the three dimensional stochastic Navier--Stokes  and Euler equations, the available results are very limited. Prior to the present work, the state of the art could be described as follows:
\begin{itemize}
\item Due to the lack of dissipation, the mere  existence of statistically stationary solutions to the three dimensional stochastic Euler equations is completely open. The subsequent and significantly more delicate step, namely the construction of ergodic stationary solutions and the question of their (non)uniqueness, seems far out of reach.
\item Nothing is known about the vanishing viscosity limit in the framework of stationary solutions in three dimensions.
\item While the existence of stationary solutions to the stochastic Navier--Stokes equations is well-established (see \cite{FG95}), their (non)uniqueness remains an outstanding open problem. The only available result in this  direction is the unique ergodicity for every Markov selection, as shown in  \cite{DD03} and \cite{FR08}. However, the uniqueness in these works only relates to the Markov process constructed therein, but as shown in \cite{HZZ21markov}, there are other Markov processes associated to the same equation  possibly having different invariant measures.
\end{itemize}

\subsection{Main results}
In this work, we are, for the first time, able to address the above problems (ii) and (iii) as well as the non-uniqueness part in point (i), in the physically relevant context of the stochastic Navier--Stokes and Euler equations  on the three dimensional torus $\mathbb{T}^3$, driven by an additive stochastic noise. The stochastic Navier--Stokes equations  read as
\begin{equation}
\label{1}
\aligned
 \dif u+\div(u\otimes u)\,\dif t+\nabla P\,\dif t&=\nu \Delta u \,\dif t+\dif B,
\\
\div u&=0,
\endaligned
\end{equation}
whereas the Euler equations are
\begin{equation}
\label{2}
\aligned
 \dif u+\div(u\otimes u)\,\dif t+\nabla P\,\dif t&=\dif B,
\\
\div u&=0.
\endaligned
\end{equation}
In the above,  $P$ is the associated pressure, $\nu>0$ is the viscosity, $B$ is a $GG^*$-Wiener process on some probability space $(\Omega, \mathcal{F}, \mathbf{P})$ and $G$ is a Hilbert--Schmidt operator from $U$ to $L^{2}_{\sigma}$ for some Hilbert space $U$ and $L^{2}_{\sigma}$  the subspace of $L^{2}$ containing mean and divergence free functions.

The primary challenge in establishing the vanishing viscosity limit and the existence of stationary solutions to the Euler equations lies in obtaining uniform in $\nu$ regularity estimates that remain valid globally in time. Traditional methods do not readily help with this task. In this paper, we propose a novel approach based  on the method of convex integration in stochastic setting, combined with the Krylov--Bogoliubov argument applied on the path space. Convex integration is an iterative procedure that constructs (evolutionary) solutions scale by scale, independently of the dissipation given by the Laplacian. Our stochastic variant represents the first instance where convex integration is successfully employed to construct statistically stationary solutions to the stochastic Euler equations.

		 In the past, convex integration has yielded numerous breakthroughs in the deterministic setting concerning the Navier--Stokes and Euler equations, as evidenced by several notable works (see e.g. \cite{BDSV,BCV18,BDLIS16, BMS20, BV19a,  CL20, DelSze2, DelSze3, DelSze13, DS17,  Ise,  L19, NV23, GKN23, GR23}). Unlike our previous efforts utilizing convex integration for these equations in a stochastic framework \cite{HZZ19, HZZ20, HZZ21markov, HZZ21},  we are now inspired by \cite{CDZ22} and overcome limitations associated with stopping times  previously used to control noise terms in the iteration.
Specifically, we eliminate the use of stopping times and instead incorporate expectations into the iterative estimates. While this approach is intuitively appealing for bypassing the challenges posed by stopping times, its implementation is far from straightforward. This difficulty arises from the quadratic nonlinearity inherent in the equations, which leads to superlinear estimates requiring meticulous analysis for each bound.

We emphasize  that the argument in \cite{CDZ22} does not inherently provide any
global-in-time $H^{\vartheta}$-bounds for some $\vartheta>0$ uniformly along the vanishing viscosity limit, which are crucial for constructing stationary solutions to the stochastic Euler equations. In our current work, we have achieved these critical estimates by carefully examining  how each $m$th moment at  the $q$th step of the iteration depends on $m$ and on the stochastic terms. This thorough analysis ultimately yields the desired global-in-time $H^{\vartheta}$-bound that remains uniform in the viscosity $\nu$.

 As we are interested in the long time behavior of solutions and particularly in the construction of statistically stationary solutions, we work with entire solutions solving the equations  for all  times $t\in\mR$. Accordingly, the norms in the convex integration scheme must be  chosen appropriately in order to provide the necessary global-in-time estimates. This is achieved through bounds of the form
$$
\sup_{\nu\geq0}\sup_{t\in\R}\mathbf{E}\left[\sup_{t\leq s \leq t+1}\|u^{\nu}(s)\|_{H^{\vartheta}}^{2r}\right]<\infty,
$$
with $r>1$ and $\vartheta>0$.
Such bounds provide uniform moment estimates locally in $C(\R;H^{\vartheta})$ and ensure the convergence of the corresponding ergodic averages,  even in the context of the stochastic Euler equations. This leads to the existence of stationary solutions.

Within this study, we focus on analytically weak solutions that satisfy the equations in the following sense.

\begin{definition}\label{d:sol}
We say that $((\Omega,\mathcal{F},(\mathcal{F}_{t})_{t\in\mR},\mathbf{P}),u,B)$ is an analytically weak solution to the Navier--Stokes system \eqref{1} provided
\begin{enumerate}
\item $(\Omega,\mathcal{F},(\mathcal{F}_{t})_{t\in\mR},\mathbf{P})$ is a stochastic basis with a complete right-continuous filtration;
\item $B$ is an $\R^{3}$-valued, spatial mean and divergence free, two-sided trace-class Brownian motion with respect to the  filtration $(\mathcal{F}_{t})_{t\in\mR}$;
\item the velocity $u \in L^2_{\rm{loc}}(\mR;L^2_{\sigma})\cap C(\mR;H^{-\delta})$ $\mathbf{P}$-a.s. for some $\delta>0$ and  is an $(\mathcal{F}_{t})_{t\in\mR}$-adapted;
\item for every $-\infty<s\leq t<\infty$ it holds  $\mathbf{P}$-a.s.
$$
\begin{aligned}
&\langle u(t),\psi \rangle + \int_{s}^{t}\langle \div(u\otimes u),\psi \rangle \dif r\ =\langle u(s),\psi \rangle +\nu\int_{s}^{t} \langle \Delta u, \psi\rangle \,\dif r +\langle B(t)-B(s),\psi\rangle
\end{aligned}
$$
for all $\psi\in C^{\infty}(\mathbb{T}^{3}),$ $\div\psi=0$.
\end{enumerate}
\end{definition}

Analytically weak solutions to the Euler equations \eqref{2} are defined in exactly the same manner, utilizing  the same function spaces, with the sole distinction  that $\nu=0$.
We observe that solutions are not required  to belong to $L^{2}_{\rm{loc}}(\mR,H^{1})$.
If $\nu>0$, they do not satisfy the corresponding energy inequality, obtained formally by testing the equations with the solution itself. In other words, the solutions to the Navier--Stokes equations we obtain are generally not the so-called Leray solutions, an issue inherent to all convex integration schemes for the Navier--Stokes equations.
However, given that the solutions constructed in this work have a prescribed  expectation of  kinetic energy, by choosing it nonincreasing they satisfy the energy inequality in the case of the Euler equations.

Given the established pathwise non-uniqueness, non-uniqueness in law, and even non-uniqueness of Markov selections in our previous works \cite{HZZ19,HZZ20, HZZ21markov}, we interpret stationarity as  shift invariance of the laws of solutions on the space of trajectories, as also discussed \cite{BFHM19, BFH20e,FFH21,HZZ22}. Specifically, we define the joint trajectory space for the solution and the driving Brownian motion as
$$
\cT = C(\mR;L^2_{\sigma})\times C(\mR;L^2_{\sigma}).
$$
This shifts $S_t$, $t\in\mR$,  on trajectories are defined by
$$
S_t(u,B)(\cdot)=(u(\cdot+t),B(\cdot+t)-B(t)),\quad t\in\mR,\quad (u,B)\in\cT.
$$
It is important to  note that the shift in the second component is adjusted to ensure that if $B$ is a Brownian motion, then $S_{t}B$ remains a Brownian motion. This construction ensures the proper shift invariance necessary for studying the stationarity of laws of solutions with respect to the underlying Brownian motion.

Stationary solutions to the stochastic Navier--Stokes equations \eqref{1} are defined as follows.

\begin{definition}\label{d:1.1}
	We say that  $((\Omega,\mathcal{F},(\mathcal{F}_{t})_{t\in\mR},\mathbf{P}),u,B)$ is a stationary solution to the stochastic Navier--Stokes equations \eqref{1} provided it satisfies  \eqref{1} in the sense of Definition~\ref{d:sol}  and its law is shift invariant, that is,
	$$\mathcal{L}[S_{t}(u,B)]=\mathcal{L}[u,B]\qquad\text{ for all }\quad t\in\mR.$$
\end{definition}

We shall distinguish this setting from the conventional invariance with respect to a Markov semigroup. The latter typically applies in scenarios where uniqueness holds, implying the Markov property is satisfied. In such cases, constructing invariant measures additionally requires the Feller property, ensuring continuous dependence on initial conditions. However, for the Navier--Stokes and Euler equations where non-uniqueness and even non-uniqueness of Markov selections have been established, we adopt this broader notion of invariance with respect to shifts on trajectories. This approach additionally allows us to bypass the strict requirements of the Feller property. Indeed, 
unlike the Feller property, the continuity of shift operators $S_{t}$ is straightforwardly ensured by the structure of the trajectory spaces $\cT$. Therefore, the use of shifts on trajectories offers a more flexible and robust approach to studying the stationarity of solutions to the Navier--Stokes and Euler equations, accommodating their pathwise non-uniqueness and avoiding unnecessary constraints associated with the Feller property.

Every stationary solution $(u,B)$ defines a dynamical system
$(\mathcal{T}, \mathcal{B} (\mathcal{T}), (S_t, t \in \mR), \mathcal{L} [u,B])
$ in the sense of, for instance, 
 \cite[Chapter 1]{DPZ96}, where $\mathcal{B}(\cT)$ denotes the $\sigma$-algebra of Borel sets in $\cT$. Consequently,  ergodicity of stationary solutions can be formulated as the ergodicity of the associated dynamical system. Furthermore,  with this notion of invariance, the existence of an ergodic stationary solution, as defined below, implies the validity of the so-called ergodic hypothesis. This hypothesis asserts that ergodic averages along trajectories of the ergodic solution converge to the ensemble average given by its law. This consideration leads us to the following definition.

\begin{definition}\label{d:1.2}
A stationary solution $((\Omega,\mathcal{F},(\mathcal{F}_{t})_{t\in\mR},\mathbf{P}),u,B)$ is ergodic provided
$$ \mathcal{L}[u,B] (A) = 1 \quad \text{or} \quad \mathcal{L}[u,B] (A) = 0 \quad \text{for all } A
\subset\mathcal{T} \text{ Borel and shift invariant}. $$
\end{definition}

The same definitions  also apply in the context of the stochastic Euler equations \eqref{2}. Specifically, the trajectory space and  regularity of the solutions are identical to those of  the Navier--Stokes equations.

With these definitions established, we can present  our first main result, addressing issues (i) and (iii) above. Rather surprisingly, the result is independent of the value of the viscosity $\nu$. Notably, it holds uniformly in the limit as viscosity tends to zero   $\nu\to 0$, encompassing both the  stochastic Navier--Stokes \eqref{1} and the Euler equations \eqref{2}. However,  in the case of the Euler equations we require higher regularity of the Brownian motion, specifically assuming $\tr((-\Delta)^{\sigma}GG^{*})<\infty$ for some $\sigma>0$.
The precise statement of the  result is provided and proved  in  Theorem~\ref{th:s1}, Theorem~\ref{thm:4.5} and Theorem \ref{th:5.1}.

\bt\label{th:main}
There exist
\begin{enumerate}
\item infinitely many stationary solutions;
\item  infinitely many ergodic stationary solutions;
\end{enumerate}
to the stochastic Navier--Stokes  \eqref{1} and Euler \eqref{2} equations. Moreover, the solutions belong to $C(\mR,H^{\vartheta})\cap C^{\vartheta}(\mR,L^{2})$ a.s. for some $\vartheta>0$.
\et

Regarding the vanishing viscosity limit (ii) as formulated above, we are able to make the following conclusion which  is proved in Theorem \ref{th:5.1} and Theorem \ref{th:new1}.

\bt\label{th:main2}
Assume $\tr((-\Delta)^{\sigma}GG^{*})<\infty$ for some $\sigma>0$. There exists $K_{0}>0$ so that for every $K\geq K_{0}$ the following holds: For an arbitrary sequence of vanishing viscosities $\nu_n\to0$, $n\in\N$, there exist a sequence of stationary solutions $u_n$, $n\in\N$, to the following stochastic Navier--Stokes equations
\begin{align*}
	\dif u_n+\div(u_n\otimes u_n)\,\dif t+\nabla P_{n}\,\dif t=\nu_n\Delta u_n\,\dif t+\dif B,
\end{align*}
so that the corresponding family of laws $\mathcal{L}[u_{n}]$, $n\in\N$, is tight in $C(\mR;L^2_\sigma)$  and
every accumulation point is a stationary solution to the stochastic Euler equations \eqref{2} satisfying
$$
\E\|u\|_{L^{2}}^{2}=K.
$$
Furthermore, the law of every stationary solution to the stochastic Euler equations \eqref{2} with some finite moment locally in $C(\mR,H^{\vartheta})\cap C^{\vartheta}(\mR,L^{2})$ for some $\vartheta>0$
 is a limit of the law of stationary solutions to the stochastic Navier--Stokes equations \eqref{1} with vanishing viscosities.
\et

\bigskip

If  $G=0$, the above results Theorem~\ref{th:main}, Theorem~\ref{th:main2}   extend to the case of deterministic Navier--Stokes \eqref{eq:detns} and Euler \eqref{eq:deteu} equations with zero external force. Moreover, with a slight adjustment  in the proof, additional insights can be derived. These 
modifications and extensions are elaborated in the proofs of Theorem \ref{th:6.21} and Theorem~\ref{thm:6.4}. The latter result also includes further observations concerning the behavior as viscosity tends to zero, noise diminishes, or both conditions occur simultaneously.

\bt\label{th:main4}
Let $G=0$.  Let  $\varepsilon>0$, $r>1$ be given and let  $ Z$ be  a stationary stochastic process with smooth trajectories, vanishing mean and divergence and satisfying
		\begin{equation*}
	\left({\mathbf{E}}\| Z\|^m_{L^2}+\E\| Z\|^m_{C^2_{t,x}}\right)^{1/m}\leq m^{1/2}L,
\end{equation*}
for any $m > 1$ and some $L \geq (2\pi)^3$. Then up to a change of probability space,
\begin{equation}\label{eq:999}
			\mathbf{E}\left[\|u-Z\|^{r}_{{C_{t}W^{1,1}}}\right]\leq\varepsilon.
		\end{equation}
		holds true for
\begin{enumerate}
\item[(1)]  the  stationary as well as ergodic stationary solutions $u$ obtained in Theorem~\ref{th:main};
		\item[(2)] the limit stationary solutions $u$ to the Euler equations \eqref{eq:deteu} obtained in Theorem~\ref{th:main2}.
		\end{enumerate}
		In particular, the solutions can be random and time dependent.
		If $Z$ is uniformly bounded in $\omega$ in $C^2_{t,x}$  then \eqref{eq:999} holds pathwise, not only in expectation.

		\et

In the proof of the above result  we make use of our stochastic convex integration construction. The key contribution  lies particularly in the  assertion \eqref{eq:999}. It shows that the solutions can exhibit certain statistics that are close to those of any prescribed process $Z$. Notably, $Z$ can be chosen to be  Gaussian or non-Gaussian.
If we consider the case of deterministic Euler equations and omit this requirement, no new convex integration construction would be necessary. Instead, we can use existing  explicit smooth steady state solutions (i.e. time independent), together with a Krein--Milman argument as in our proof  to deduce the result of Theorem~\ref{th:main}  in this Euler case with $G=0$.
On the other hand, applying the same strategy to  the deterministic  Navier--Stokes equations is  more delicate is more intricate due to dissipation effects. In four  dimensions, the steady state solutions identified in  \cite{L19} can be employed similarly to establish the
 existence of non-unique  ergodic statistically stationary solutions.

\subsection*{Organization of the paper} In Section~\ref{s:not}, we collect the basic notations used throughout the paper. Section~\ref{s:in} is the core of our proofs: here, the stochastic convex integration is developed and employed to construct entire non-unique analytically weak and probabilistically strong solutions with a prescribed kinetic energy. This is then used in Section~\ref{s:4} together with a Krylov--Bogoliubov's argument to obtain existence of non-unique stationary solutions to the stochastic Navier--Stokes equations. The results concerning stationary solutions to the stochastic Euler equations, the vanishing viscosity limit can be found in Section~\ref{s:5}, whereas the results for the deterministic systems are proved in Section~\ref{s:6}. In Appendix~\ref{s:B}, we recall the construction of intermittent jets from \cite{BCV18, BV19} and in Appendix \ref{ap:B} we give estimates on  amplitude functions used in the convex integration construction.

\section{Notations}
\label{s:not}

\subsection{Function spaces}
  Throughout the paper, we employ the notation $a\lesssim b$ if there exists a constant $c>0$ such that $a\leq cb$, and we write $a\simeq b$ if $a\lesssim b$ and $b\lesssim a$. We let $\mN_{0}:=\mN\cup \{0\}$. Given a Banach space $E$ with a norm $\|\cdot\|_E$ and $t\in\mR$, we write $C_tE=C([t,t+1];E)$ for the space of continuous functions from $[t,t+1]$ to $E$, equipped with the supremum norm $\|f\|_{C_tE}=\sup_{s\in[t,t+1]}\|f(s)\|_{E}$.  For $\alpha\in(0,1)$ we  define $C^\alpha_tE$ as the space of $\alpha$-H\"{o}lder continuous functions from $[t,t+1]$ to $E$, endowed with the norm $\|f\|_{C^\alpha_tE}=\sup_{s,r\in[t,t+1],s\neq r}\frac{\|f(s)-f(r)\|_E}{|r-s|^\alpha}+\sup_{s\in[t,t+1]}\|f(s)\|_{E}.$ Here we use $C_t^\alpha$ to denote the case when $E=\mathbb{R}$. {We also write $C_b(\mR;E)$ for functions in $C(\mR;E)$ such that $\|f\|_{C_{b}(\mR;E)}:=\sup_{t\in \mR}\|f(t)\|_E<\infty$. For $\beta\in (0,1]$ we define $C_b^\beta(\mR;E)$ as functions in $C^\beta(\mR;E)$ such that $\|f\|_{C_b^\beta(\mR;E)}:=\sup_{t\in \mR}\|f\|_{C^\beta_tE}<\infty$.}
    We use $L^p$ to denote the set of  standard $L^p$-integrable functions from $\mathbb{T}^3$ to $\mathbb{R}^3$. For $s>0$, $p>1$ we set $W^{s,p}:=\{f\in L^p; \|(I-\Delta)^{{s}/{2}}f\|_{L^p}<\infty\}$ with the norm  $\|f\|_{W^{s,p}}=\|(I-\Delta)^{{s}/{2}}f\|_{L^p}$. Set $L^{2}_{\sigma}=\{f\in L^2; \int_{\mathbb{T}^{3}} f\,\dif x=0,\div f=0\}$. For $s>0$, we define $H^s:=W^{s,2}\cap L^2_\sigma$. For $s<0$ we define $H^s$ to be the dual space of $H^{-s}$.
For $t\in\mathbb{R}$ and a domain $D\subset\R^{+}$ we denote by  $C^{N}_{t,x}$ and $C^{N}_{D,x}$, respectively, the space of $C^{N}$-functions on $[t,t+1]\times\mathbb{T}^{3}$ and on $D\times\mathbb{T}^{3}$, respectively,  $N\in\N_{0}$. The spaces are equipped with the norms
$$
\|f\|_{C^N_{t,x}}=\sum_{\substack{0\leq n+|\alpha|\leq N\\ n\in\N_{0},\alpha\in\N^{3}_{0} }}\|\partial_t^n D^\alpha f\|_{L^\infty_{[t,t+1]} L^\infty},\qquad \|f\|_{C^N_{D,x}}=\sum_{\substack{0\leq n+|\alpha|\leq N\\ n\in\N_{0},\alpha\in\N^{3}_{0} }}\sup_{t\in D}\|\partial_t^n D^\alpha f\|_{ L^\infty}.
$$
 For a Polish space $H$ we denote by $\mathcal{B}(H)$  the $\sigma$-algebra of Borel sets in $H$.  We also use $\mathring{\otimes}$ to denote the trace-free part of the tensor product.

By $\mathbb{P}$ we denote the Helmholtz projection. We recall the inverse divergence operator $\mathcal{R}$ from \cite[Section 5.6]{BV19}, which acts on vector fields $v$ with $\int_{\mathbb{T}^3}v\dif x=0$ as
\begin{equation*}
	(\mathcal{R}v)^{kl}=(\partial_k\Delta^{-1}v^l+\partial_l\Delta^{-1}v^k)-\frac{1}{2}(\delta_{kl}+\partial_k\partial_l\Delta^{-1})\div\Delta^{-1}v,
\end{equation*}
for $k,l\in\{1,2,3\}$. Then $\mathcal{R}v(x)$ is a symmetric trace-free matrix for each $x\in\mathbb{T}^3$, and $\mathcal{R}$ is a right inverse of the div operator, i.e. $\div(\mathcal{R} v)=v$. By \cite[Theorem B.3]{CL20} we know
\begin{align}\label{eR}\|\mathcal{R}f(\sigma\cdot)\|_{L^p}\lesssim \sigma^{-1}\|f\|_{L^p}\quad\textrm{for } \sigma\in\mathbb{N}.\end{align}

By $\cS^{3\times 3}$ we denote the set of symmetric $3\times 3$ matrices and by $\cS_0^{3\times 3}$  the set of  symmetric trace-free matrices. Let  $C_0^\infty(\mT^3,\mathbb R^{3\times 3})$ be the set of periodic smooth matrix valued functions with zero mean. We also introduce the bilinear version $\cB:C^\infty(\mT^3,\mR^3)\times C_0^\infty(\mT^3,\mathbb R^{3\times 3})\to C^\infty(\mT^3,\cS_0^{3\times 3})$ as in \cite[Section~B.3]{CL20} by
\begin{align*}
	\cB(v,A)=v\cR A-\cR(\nabla v\cR A).
\end{align*}
Then by \cite[Theorem B.4]{CL20} we have $\div(\mathcal{B} (v,A))=vA-\frac1{(2\pi)^3}\int_{\mT^3} vA\dif x $ and
\begin{align}\label{eB}
	\|\cB(v,A)\|_{L^p}\lesssim \|v\|_{C^1}\|\mathcal{R}A\|_{L^p}.
\end{align}

\subsection{Probabilistic elements}\label{s:2.2}

 For a given probability measure $P$ we denote by $\mathbf{E}^P$  the expectation under $P$.
Regarding the driving noise, we assume that $B$ is an $\R^{3}$-valued two-sided  $GG^*$-Wiener process with zero spatial mean and zero divergence, defined on some probability space $(\Omega, \mathcal{F}, \mathbf{P})$ and $G$ is a Hilbert--Schmidt operator from $U$ to $L_{\sigma}^2$ for some Hilbert space $U$.

For $p\in[1,\infty)$ we denote
	$$
	\$u\$^p_{L^2,p}:=\sup_{t\in \mR}\mathbf{E}\left[\sup_{t\leq s\leq t+1}\|u(s)\|^p_{L^2}\right],\quad \$u\$^p_{C^1_{t,x},p}:=\sup_{t\in \mR}\mathbf{E}\left[\|u(s)\|^p_{C^1_{[t,t+1],x}}\right].
	$$
	These norms define function spaces of random variables on $\Omega$ taking values in $C(\mR,L^{2})$ and {$C^{1}(\mR\times\mathbb{T}^{3})$}, respectively, with  bounds in $L^{p}(\Omega;C(I,L^{2}))$ and $L^{p}(\Omega;{C^{1}(I\times\mathbb{T}^{3}))}$ for any bounded interval $I\subset\mR$. Furthermore, the bounds  only depend on the length of the interval $I$, not on its location within $\mR$. In the sequel, we simply say that $u$ has a uniform moment of order $p$ locally in $C(\mR;L^{2})$ provided $\$ u\$_{L^{2},p}<\infty$.
	Similarly, we define the corresponding norms with $L^2$ replaced by $L^p$, $H^\vartheta$ or $C^1_{t,x}$ replaced by  $C_t^{\frac{1}{2}-2\delta}L^\infty$, $C_t^\vartheta L^2$ and $C_tW^{1,p}$.

\section{Stochastic convex integration}
\label{s:in}

The previous works using convex integration in the stochastic setting always reduced the problem to the deterministic setting by introducing suitable stopping times.\footnote{The first exception was our previous work on a class of  supercritical/critical SPDEs with an irregular spatial perturbation \cite{HZZ22}. Due to the time independence of the noise, no stopping times were necessary.}
This permitted to control the noise uniformly in $\omega$ so that the convex integration could proceed pathwise up to the stopping time. The stopping times were then removed a posteriori by a suitable extension of solutions. Such an approach is not suitable for the construction of stationary solutions. Hence, inspired by \cite{CDZ22} we present an honest stochastic convex integration, constructing directly solutions on the whole time line $\R$. This is achieved by introducing expectations to the iterative estimates in convex integration. The main difficulty lies in the fact that due to the quadratic nonlinearity the estimates are superlinear. More precisely, the estimate of any $p$th moment at the level $q+1$  necessarily contains higher moments at the level $q$. Accordingly, all the estimates need to be tracked down very carefully, paying a particular attention to the appearing  constants, what they depend on and how  precisely. Otherwise, it would not be possible to close the estimates. The key observation is that the superlinear terms always contain a small constant which may be used to absorb the bounds.

As the first step, we decompose a solution to the Navier--Stokes system \eqref{1} with $\nu=1$ into the sum $u=z+v$ where $z$ is the unique stationary solution to the linear stochastic heat equation
\begin{equation}\label{linear}
\aligned
\dif z -(\Delta-1) z \,\dif t&=\dif B,
\endaligned
\end{equation}
where $B$ is a $\mathbb{R}^3$-valued two-sided trace-class Wiener process with spatial zero mean  (see e.g. \cite[page 99]{PR07}),
and $v$ solves the nonlinear deterministic equation
\begin{equation}\label{nonlinear}
\aligned
\partial_tv -\Delta v-z+\div((v+z)\otimes (v+z))+\nabla P&=0,
\\
\div v&=0.
\endaligned
\end{equation}
Here, $z$ is divergence free by the assumptions on the noise and by $P$ we denote the  pressure term associated to $v$.

\begin{remark}
For notational simplicity, we work in this section as well as in Section~\ref{s:4}  with the unit viscosity $\nu=1$. This fact is used only in Proposition~\ref{fe z} below, which profits from the smoothing effect of the Laplacian. More precisely, the spatial regularity is needed for the convergence rate in the convex integration in order to deduce the  $H^\vartheta$-estimate. Only a bound in $L^2$ would not be enough. For a general $\nu$, the bound in Proposition~\ref{fe z} would depend on $\nu^{-p}$. Hence, for the results regarding stationary solutions to the Euler equations,  the vanishing viscosity limit  and the anomalous dissipation in Section~\ref{s:5} and Section~\ref{s:6}, it is necessary to increase the regularity of the noise in order to compensate for the lack of smoothing of the linear part.
\end{remark}

We included the linear zero order term $z$ on the left hand side of \eqref{linear} in order to obtain a stabilization of the equation needed for the necessary global in time estimates. It will be seen in the course of the construction that the corresponding counter term in \eqref{nonlinear} will not cause any difficulties.
Using the  factorization method it is standard to derive  regularity of the stochastic convolution $z$  on a given stochastic basis $(\Omega, \mathcal{F},(\mathcal{F}_{t})_{t\in\mR},\mathbf{P})$ with $(\mathcal{F}_{t})_{t\in\mR}$ canonical filtration given in \cite[page 99]{PR07}.  In particular, the following result follows from \cite[Theorem 5.16]{DPZ92} together with the Kolmogorov continuity criterion.

\bp\label{fe z}
	Suppose that $\mathrm{Tr}(GG^*)<\infty$. Then for any $\delta\in (0,1/2)$, $p\geq2$
\begin{equation}\label{eq:u0}
	\sup_{t\in\mR}\mathbf{E}\left[\sup_{t\leq s\leq t+1}\|z(s)\|^p_{H^{1-\delta}}+\|z\|^p_{C_{t}^{1/2-\delta}L^2}\right]\leq (p-1)^{p/2}L^p,
\end{equation}
where $L\geq 1$ depends  on $\mathrm{Tr}(GG^{*})$, $\delta$ and is independent of $p$.
	\ep

\begin{proof}
We recall that the unique stationary solution to \eqref{linear} has the explicit form $z(t)=\int_{-\infty}^t S(t-s)\dif B(s)$ where $S(t)=e^{t(\Delta-I)}$, $t\geq0$. The Wiener process $B$ is given by $B=\sum_{k\in\mN}c_ke_k\beta_k$ for an orthonormal basis   $\{e_k\}_{k\in\mN}$ of $L^{2}_{\sigma}$, a sequence of mutually independent standard two-sided Brownian motions $\{\beta_{k}\}_{k\in\mN}$ and the coefficients satisfy $\sum_{k\in\mN}c_k^2<\infty$. Then it holds for $\gamma\in (0,1/2)$, $t\geq s$
	\begin{align*}
		\E\|z(t)-z(s)\|_{L^2}^2&=\sum_{k=1}^{\infty}c_k^2\int_s^t\|S(t-\sigma)e_k\|_{L^2}^2\dif \sigma+\sum_{k=1}^\infty\int_{-\infty}^sc_k^2\|[S(t-\sigma)-S(s-\sigma)]e_k\|^2_{L^2}\dif \sigma
		\\&\leq M\tr(GG^*)\Big[(t-s)+\int_{-\infty}^se^{-2(s-\sigma)}\Big|\int_{s-\sigma}^{t-\sigma}\frac{1}r\dif r\Big|^2\dif \sigma\Big]
		\\&\leq M\tr(GG^*)\Big[(t-s)+\int_{-\infty}^se^{-2(s-\sigma)}(s-\sigma)^{-2\gamma}\Big|\int_{s-\sigma}^{t-\sigma}r^{\gamma-1}\dif r\Big|^2\dif \sigma\Big]
		\\&\leq M\tr(GG^*)[(t-s)+(t-s)^{2\gamma}],
	\end{align*}
	where the constant $M$  depends only on the semigroup and $\gamma$ but is independent of time.
Using Gaussianity we have
\begin{align*}
\E\|z(t)-z(s)\|_{L^2}^{p}\leq (p-1)^{p/2}\Big(\E\|z(t)-z(s)\|_{L^2}^{2}\Big)^{p/2}.
\end{align*}
Similar computations can  be performed for the $H^{1-\delta}$-norm as well, only the resulting time regularity is lower and depends on $\delta$.
	The  result then follows from Kolmogorov's continuity criterion.
\end{proof}

In the following  we choose $L\geq (2\pi)^{3/2}$ for simplicity.

Global in time estimates of the form \eqref{eq:u0} are well-suited for the application of a Krylov--Bogoliubov argument leading to existence of stationary solutions as limits of ergodic averages. Our goal in this section is to construct solutions to the Navier--Stokes system \eqref{1} satisfying similar bounds. To this end,
we use the norms $\$\cdot\$_{L^2,p}$ and $\$\cdot\$_{C^1_{t,x},p}$ introduced in Section \ref{s:2.2}, which play the essential role in the construction.

Let us now explain how the convex integration iteration is set up.
We consider an increasing sequence $\{\lambda_q\}_{q\in\mathbb{N}_{0}}\subset \mathbb{N}$ which diverges to $\infty$, and a sequence $\{\delta_q\}_{q\in \mathbb{N}}\subset(0,1)$  which is decreasing to $0$. We choose $a\in\mathbb{N},$ $ b\in\mathbb{N},$ $  \beta\in (0,1]$ and let
$$\lambda_q=a^{(b^q)}, \quad\delta_1=1,\quad \delta_q=\frac12\lambda_2^{2\beta}\lambda_q^{-2\beta},\quad q\geq2.$$
Here $\beta$ will be chosen sufficiently small and $a$ as well as $b$ will be chosen sufficiently large.
We assume $\sum_{r\geq 1}\delta_{r}^{1/2}\leq 1+\sum_{r\geq2}a^{{b^{2}}\beta-(r-1)b^2\beta}=1+\frac{1}{1-a^{-\beta b^2}}\leq 3$ which boils down  to
\begin{equation}\label{aaa}
{a^{\beta b^2}\geq 2}.
\end{equation}
Here we require more than is necessary for the later use and we keep this assumption  from now on.
 More details on the choice of these parameters will be given below in the course of the construction.
The iteration is indexed by a parameter $q\in\mathbb{N}_{0}$. At each step $q$, a pair $(v_q, \mathring{R}_q)$ is constructed solving the following system
\begin{equation}\label{induction ps}
\aligned
\partial_tv_q-z_q-\Delta v_q +\div((v_q+{z_q})\otimes (v_q+{z_q}))+\nabla p_q&=\div \mathring{R}_q,
\\
\div v_q&=0.
\endaligned
\end{equation}
In the above we   define $z_q=\mP_{\leq f(q)}z$ with $f(q)=\lambda_{q+1}^{\alpha/8}$ and $\mathring{R}_q$ is trace-free and we put the trace part into the pressure. Here $\alpha\in(0,1)$ is a small parameter and precise definition is given in Section \ref{s:par}. Thanks to this approximation of $z$, we are able to lower the assumption on the spatial regularity of the noise $B$, namely, to cover the case of any trace-class noise.
We observe that
\begin{equation}\label{z ps}
\$ z_q\$_{L^\infty,p}\leq (p-1)^{1/2}L\lambda_{q+1}^{\alpha/8},
 \quad
 \$z_q\$_{C_t^{\frac{1}{2}-2\delta}L^\infty,p}\leq (p-1)^{1/2}L\lambda_{q+1}^{\alpha/4}.
\end{equation}

We intend to  construct approximations $v_{q}$ with a uniform moment of order $2r$ for a given $r>1$ locally in $C(\mR,H^{\vartheta})$ and $C^{\vartheta}(\mR,L^{2})$ for some $\vartheta>0$, in the sense of the norms $\$\cdot\$_{H^{\vartheta},2r}$ and $\$\cdot\$_{C^{\vartheta}_tL^{2},2r}$. But to this end, it is  necessary to quantify the blow up of the higher moments, as these also appear in the estimates.

Under the above assumptions, our main iteration  reads as follows, the proof of this result is presented in Section~\ref{s:it} below.

\begin{proposition}\label{p:iteration}
Assume \eqref{eq:u0} and let $r>1$ be fixed. Given smooth function $e:\R\to(0,\infty)$ so that $\bar e\geq e(t)\geq \underline{e}\geq1 $  with $\|e'\|_{C^{0}}\leq \tilde e$ for some constants $\overline{e},\underline{e},\tilde{e}>0$, there exists a choice of parameters $a, b, \beta$ and $\alpha\in(0,1/49)$ with $\alpha b>32/7$  such that the following holds true: Let $(v_{q},\mathring{R}_{q})$ for some $q\in\N_{0}$ be an $(\mathcal{F}_{t})_{t\geq 0}$-adapted solution to \eqref{induction ps} satisfying
\begin{equation}\label{inductionv}
\$v_{q}\$_{L^{2},2r}\leq
M_0 \bar e^{1/2}\sum_{ k=1}^{q}\delta_{k}^{1/2}
\end{equation}
 for a universal constant $M_0\geq 1$, and for $m\geq 1$
\begin{align}\label{inductionv m}
 \$v_q\$_{L^2,m}&\leq M_0(6^{q-1}\cdot 12mL^2)^{3(6^{q-1})}{\lambda_q}+M_0\bar e^{1/2}\sum_{ r=1}^{q}\delta_{r}^{1/2},
\end{align}
and
\begin{align}\label{inductionv C1}
 \$v_q\$_{C^1_{t,x},m}&\leq \lambda_{q}^{23/7}(6^{q-1}\cdot 16mL^2)^{4(6^{q-1})},\quad \$v_q\$_{C^2_{t,x},m}\leq \lambda_{q}^{37/7}(6^{q-1}\cdot 20mL^2)^{5(6^{q-1})},
\end{align}
\begin{align}\label{eq:R}
\$\mathring{R}_q\$_{L^1,r}\leq
\frac1{48}\delta_{q+2}\underline{e},
\end{align}
and for $m\geq1$
\begin{align}\label{bd:R}
 \$\mathring{R}_q\$_{L^1,m}\leq (6^q\cdot 4mL^2)^{(6^q)}{\lambda^2_q}.
\end{align}
 Moreover, for any $t\in\mR$
\begin{align}\label{p:gamma}
\frac34\delta_{q+1}e(t)\leq e(t)-\E\|(v_q+z_q)(t)\|_{L^2}^2\leq \frac54\delta_{q+1}e(t),
\end{align}
 Then    there exists an $(\mathcal{F}_{t})_{t\geq 0}$-adapted process $(v_{q+1},\mathring{R}_{q+1})$ which solves \eqref{induction ps}, obeys \eqref{inductionv}-- \eqref{p:gamma}  at the level $q+1$ and satisfies
 \begin{equation}\label{iteration}
\$v_{q+1}-v_q\$_{L^2,2r}\leq
M_0\bar e^{1/2}\delta_{q+1}^{1/2},
\end{equation}
and for $p=\frac{32}{32-7\alpha}$
\begin{align}\label{bdvq:w1p}
	\$v_{q+1}-v_q\$_{C_tW^{1,p},r}\leq {\lambda_{q+1}^{-\alpha/2}}.
\end{align}
\end{proposition}

Using  \eqref{inductionv C1} and the choice of the parameters in Section \ref{s:par} we have
\begin{align}\label{vqc1}
	\$v_q\$_{C^1_{t,x},2r}\leq \lambda_q^4,\quad \$v_q\$_{C^2_{t,x},r}\leq \lambda_q^6.
\end{align}

We start the iteration from $v_{0}\equiv 0$ on $\mR$. In that case, we have
$\mathring{R}_0=z_{0}\mathring\otimes z_{0}-\mathcal{R}z_{0}$
so that
\begin{align*}
\$\mathring{R}_0\$_{L^1,m}\leq \$z\$^2_{L^2,2m}+(2\pi)^{3/2}\$z\$_{L^2,m}\leq 4mL^2
\end{align*}
and \eqref{eq:R}, \eqref{bd:R}  are satisfied on the level  $q=0$, since $\delta_{2}=1/2$ and provided $8\cdot 48r L^{2}\leq \underline{e}$. Here, we used $L\geq (2\pi)^{3/2}$. For \eqref{p:gamma} we require
\begin{align*}
	\frac34e(t)\leq e(t)-\E\|z_0(t)\|_{L^2}^2\leq \frac54e(t),
\end{align*}
{which is satisfied provided $e(t)\geq 4L^{2}$.}

We deduce the following result.

\bt\label{thm:6.1}
Let $r>1$ and a smooth function $e:\mR\to(0,\infty)$ satisfying $\bar e\geq e(t)\geq \underline{e}\geq 8{\cdot 48}r L^2$ be given.
There exists an $(\mathcal{F}_t)_{t\in\mR}$-adapted process $u$ which belongs to $ C(\mR,H^{\vartheta})\cap C^{\vartheta}(\mR,L^{2})$ $\mathbf{P}$-a.s.  for some $\vartheta>0$  and is an analytically weak solution to \eqref{1} in the sense of Definition~\ref{d:sol}.
Moreover, the solution satisfies
\begin{align}\label{est:u1}\$u\$_{H^\vartheta,2r}+\$u\$_{C_{t}^\vartheta L^2,2r}<\infty,\end{align}
and for all $t\in\mR$
\begin{align}\label{eq:K1}
\begin{aligned}
\mathbf{E}\|u(t)\|_{L^2}^2=e(t).
\end{aligned}
\end{align}
 Furthermore, for every $\eps>0$
one may find solution such that $v=u-z$ satisfying
 \begin{align}\label{es:v1}
\$v\$_{C_{t}W^{1,1},r}\leq \eps.\end{align}
\et

\begin{proof}
By interpolation we deduce for $\vartheta\in (0,\frac{\beta}{4+\beta})$,
	\begin{equation}\label{eqtheta}
	\sum_{q\geq0}\$v_{q+1}-v_q\$_{H^{\vartheta},2}\lesssim \sum_{q\geq0}\$v_{q+1}-v_q\$_{L^2,2}^{1-\vartheta}\$v_{q+1}-v_q\$_{H^1,2}^{\vartheta}\lesssim \sum_{q\geq0}\delta_{q+1}^{\frac{1-\vartheta}{2}}\lambda_{q+1}^{4\vartheta}<\infty.
	\end{equation}
Similarly we could change $H^{\vartheta}$ to $C^\vartheta_tL^2$.
As a consequence, a limit $v=\lim_{q\rightarrow\infty}v_q$ exists and lies in $L^2(\Omega,C(\mR,H^{\vartheta})\cap C^\vartheta(\mR,L^2))$.  Since $v_q$ is $(\mathcal{F}_t)_{t\in\mR}$-adapted for every $q\in\mathbb{N}_{0}$, the limit
	$v$ is $(\mathcal{F}_t)_{t\in\mR}$-adapted as well.
	Furthermore, it follows from \eqref{eq:R} that $\lim_{q\rightarrow\infty}\mathring{R}_q=0$ in $L^1(\Omega,C(\mR;L^1))$ and $\lim_{q\rightarrow\infty}z_q=z$ in $L^p(\Omega,C(\mR;L^2))$ for any $p\geq1$. Thus $v$ is an analytically  weak solution to (\ref{nonlinear}).  Hence letting $u=v+z$ we obtain an $(\mathcal{F}_{t})_{t\in\mR}$-adapted analytically weak solution to \eqref{1}. Moreover, the estimate for $u$ holds.
Finally, 	\eqref{eq:K1} follows from \eqref{p:gamma}.

For the last result, we use \eqref{bdvq:w1p} and conditions on $\alpha$ to have
\begin{align*}
	\$v\$_{C_{t}W^{1,p},r}&\lesssim \sum_{q=0}^\infty \$v_{q+1}-v_q\$_{C_tW^{1,p},r}
	\leq \sum_{q=0}^\infty{\lambda_{q+1}^{-\alpha/2}}\lesssim \frac{a^{-\alpha b/2}}{1-a^{-\alpha b/2}}
	\\&=\frac{1}{a^{\alpha b/2}-1}\leq \frac1{a^{{16}/7}-1}\leq \eps,
\end{align*}
where we use $\alpha b>{32}/7$ and $-1/7+6\alpha<-\alpha/2$ and we may choose $a$ large enough such that the last inequality holds.
	\end{proof}

\subsection{Proof of Proposition~\ref{p:iteration}}
\label{s:it}

The proof proceeds in several main steps which are the same in many convex integration schemes. First of all, we start the construction by fixing the parameters in Section~\ref{s:par} and proceed with a mollification step in Section~\ref{s:p}. Section~\ref{s:313} introduces the new iteration $v_{q+1}$. This is the main part of the construction which differs in each convex integration scheme.   Here, we   construct new  amplitude functions $a_{(\xi)}$ similarly to \cite{HZZ21markov} but we replace the pathwise construction by a stochastic variant, namely, we work explicitly with  expectations of $v_q+z_q$. Section~\ref{sss:v} contains the inductive estimates of $v_{q+1}$, especially the moment estimates, whereas in Section~\ref{s:en} we show how the energy is controlled. Finally, in Section~\ref{s:def}, we define the new stress $\mathring{R}_{q+1}$ and  establish the inductive moment estimate on $\mathring{R}_{q+1}$ in Section~\ref{sss:R}.

\subsubsection{Choice of parameters}
\label{s:par}

In the sequel, additional parameters will be indispensable and their value has to be carefully chosen in order to respect all the compatibility conditions appearing in the estimations below. First, for a sufficiently small  $\alpha\in (0,1)$ to be chosen below, we let $\ell\in (0,1)$ be  a small parameter  satisfying
	\begin{equation}\label{ell}
	\ell \lambda_q^{4}\leq \lambda_{q+1}^{-\alpha},\quad \ell^{-1}\leq \lambda_{q+1}^{2\alpha},\quad\bar{e}\leq\ell^{-1}.
	\end{equation}
	In particular, we define
	\begin{equation}\label{ell1}
	\ell:=\lambda_{q+1}^{-{3\alpha}/{2}}\lambda_q^{-2}.
	\end{equation}
	In the sequel, we  use the following bounds
	$$
	\alpha b>32/7,\quad {43}\alpha<1/14,\quad \alpha>40\beta b^2,\quad {18}/b+2\beta b^2<1/14,\quad2\beta b <\frac1{7}-{132}\alpha,
	$$
	which can be obtained by choosing $\alpha$ small  such that
	$\frac1{{133}\cdot 7}>\alpha,$
	and choosing $b\in\mathbb{N}$ large enough such that
	$\alpha b>32/7$ and finally choosing $\beta$ small such that
	$\alpha>40\beta b^2$.
Hence, we shall choose rational $\alpha$ small first and $b$ large, then $\beta$  small enough. The last free parameter is $a$ which satisfies the lower bounds given through  \eqref{aaa} and the last bound in \eqref{ell}. Let $c$ satisfy $q6^q\leq c7^q$. We then choose $a\geq (192 L^2 r)^c\vee (252 L^2)^{3c}$.
	In the sequel, we increase $a$ in order to absorb various implicit and universal constants.

\subsubsection{Mollification}\label{s:p}
	
	We intend to replace $v_q$ by a mollified velocity field $v_\ell$. To this end, let $\{\phi_\varepsilon\}_{\varepsilon>0}$ be a family of standard mollifiers on $\mathbb{R}^3$, and let $\{\varphi_\varepsilon\}_{\varepsilon>0}$ be a family of  standard mollifiers with support in $(0,1)$. The one-sided mollifier here is used in order to preserve adaptedness. We define a mollification of $v_q$, $\mathring{R}_q$ and $z_q$ in space and time by convolution as follows
	$$v_\ell=(v_q*_x\phi_\ell)*_t\varphi_\ell,\qquad
	\mathring{R}_\ell=(\mathring{R}_q*_x\phi_\ell)*_t\varphi_\ell,\qquad
	z_\ell=({z_q}*_x\phi_\ell)*_t\varphi_\ell,$$
	where $\phi_\ell=\frac{1}{\ell^3}\phi(\frac{\cdot}{\ell})$ and $\varphi_\ell=\frac{1}{\ell}\varphi(\frac{\cdot}{\ell})$.
	Since the mollifier $\varphi_\ell$ is supported on $(0,1)$, it is easy to see that $z_\ell$ is $(\mathcal{F}_t)_{t\in\mR}$-adapted and so are $v_\ell$ and $\mathring{R}_\ell$.
	Since (\ref{induction ps}) holds, it follows that $(v_\ell,\mathring{R}_\ell)$ satisfies
	\begin{equation}\label{mollification}
	\aligned
	\partial_tv_\ell -z_\ell-\Delta v_\ell+\div((v_\ell+z_\ell)\otimes (v_\ell+z_\ell))+\nabla p_\ell&=\div (\mathring{R}_\ell+R_{\textrm{com}})
	\\\div v_\ell&=0,
	\endaligned
	\end{equation}
	where
	\begin{equation*}
	R_{\textrm{com}}=(v_\ell+z_\ell)\mathring{\otimes}(v_\ell+z_\ell)-((v_q+{z_q})\mathring{\otimes}(v_q+{z_q}))*_x\phi_\ell*_t\varphi_\ell,
	\end{equation*}
	\begin{equation*}
	p_\ell=(p_q*_x\phi_\ell)*_t\varphi_\ell-\frac{1}{3}\big(|v_\ell+z_\ell|^2-(|v_q+{z_q}|^2*_x\phi_\ell)*_t\varphi_\ell\big).
	\end{equation*}

We have
\begin{equation}\label{error}
\|v_q(t)-v_\ell(t)\|_{L^2}\lesssim\ell \|v_q\|_{C^1_{[t-1,t+1],x}},
\end{equation}
which by \eqref{vqc1} implies that
\begin{equation}\label{error ps}
\$v_q-v_\ell\$_{L^2,2r}\lesssim\ell \$v_q\$_{C^1_{t,x},2r}\leq \ell\lambda_q^{4}\leq \frac{1}{4} \bar e^{1/2}\delta_{q+1}^{1/2},
\end{equation}
where we  used the fact that $\ell\lambda_q^{4}<\lambda_{q+1}^{-\beta}$.
In addition,
\begin{equation}\label{eq:vl}
\|v_\ell\|_{C_tL^2}\leq \|v_q\|_{C_{[t-1,t+1]}L^2}.
\end{equation}

\subsubsection{Construction of $v_{q+1}$}\label{s:313}

Let us now proceed with the construction of the perturbation $w_{q+1}$ which then defines the next iteration by $v_{q+1}:=v_{\ell}+w_{q+1}$.
To this end, we employ  the intermittent jets introduced in \cite{BCV18} and presented in \cite[Section 7.4]{BV19}, which we recall in Appendix~\ref{s:B}. In particular, the building blocks $W_{(\xi)}=W_{\xi,r_\perp,r_\|,\lambda,\mu}$ for $\xi\in\Lambda$ are defined in (\ref{intermittent}) and the set $\Lambda$ is introduced in Lemma \ref{geometric}.
The necessary estimates are collected  in \eqref{bounds}.  We choose the following parameters
\begin{equation}\label{parameter}
	\aligned
	\lambda&=\lambda_{q+1},
	\qquad
	r_\|=\lambda_{q+1}^{-4/7},
	\qquad r_\perp=r_\|^{-1/4}\lambda_{q+1}^{-1}=\lambda_{q+1}^{-6/7},
	\qquad
	\mu=\lambda_{q+1}r_\|r_\perp^{-1}=\lambda_{q+1}^{9/7}.
	\endaligned
\end{equation}
It is required that $b$ is a multiple of $7$ to ensure that $\lambda_{q+1}r_\perp= a^{(b^{q+1})/7}\in\mathbb{N}$.

Now we define $\rho$ as follows
\begin{equation}\label{eq:rho3}
	\rho:=2\sqrt{\ell^2+|\mathring{R}_\ell|^2}+\gamma_\ell,
\end{equation}
$$\gamma_q(t):=\frac1{3\cdot (2\pi)^3} \Big[ e(t)(1-\delta_{q+2})-\E\|v_q(t)+{z_q}(t)\|_{L^2}^2\Big],$$
and
$$ \gamma_\ell:=\gamma_q*_t\varphi_\ell.$$
We observe that  \eqref{aaa} and $b\geq 2$ implies in particular
$
\frac43\leq a^{2\beta b(b-1)},
$
i.e. $\frac34\delta_{q+1}\geq \delta_{q+2}$, it follows that $\gamma_{q}\geq 0$.
In view of the definition of $\rho $ in \eqref{eq:rho3}, we  obtain for any $p\in [1,\infty]$,
\begin{equation}\label{rho}
	\|\rho\|_{L^p}\leq 2\ell(2\pi)^{3/p}+2\|\mathring{R}_\ell\|_{L^p}+\frac12\delta_{q+1}\bar e.
\end{equation}
Furthermore, by mollification estimates, the embedding $W^{4,1}\subset L^\infty$  we obtain for $N\geq0$
	$$
	\|\mathring{R}_\ell\|_{C^N_{t,x}}\lesssim \ell^{-4-N}\|\mathring{R}_q\|_{C_{[t-1,t+1]}L^1},
	$$
	which in particular leads to
\begin{equation}\label{rho0}
	\|\rho\|_{C^{0}_{t,x}}\lesssim \ell+\ell^{-4}\|\mathring{R}_q\|_{C_{[t-1,t+1]}L^1}+\delta_{q+1}\bar e.
\end{equation}
We put further details on the $C^N_{t,x}$-estimates of $\rho$ in Appendix~\ref{ap:B}
	and  by \eqref{rhoN1}  we obtain
\begin{align}\label{rhoN}
	\|\rho\|_{C^N_{t,x}}
	\lesssim \ell^{- 4 - N} \|\mathring{R}_q\|_{C_{[t-1,t+1]}L^1} + \ell^{- 6N + 1}
	\|\mathring{R}_q\|_{C_{[t-1,t+1]}L^1}^N+\ell^{-N}\delta_{q+1}\bar e. \end{align}

Now, we define the amplitude functions
\begin{equation}\label{amplitudes}a_{(\xi)}(\omega,t,x):=a_{\xi,q+1}(\omega,t,x):=\rho(\omega,t,x)^{1/2}\gamma_\xi\left(\Id
	-\frac{\mathring{R}_\ell(\omega,t,x)}{\rho(\omega,t,x)}\right),\end{equation}
where $\gamma_\xi$ is introduced in  Lemma \ref{geometric}.
By  (\ref{geometric equality}) we have
\begin{equation}
\label{cancellation}
(2\pi)^{-3}\sum_{\xi\in\Lambda}a_{(\xi)}^2\int_{\mathbb{T}^3}W_{(\xi)}\otimes W_{(\xi)}\dif x=\rho \Id-\mathring{R}_\ell,
\end{equation}
and using \eqref{estimate a1}
\begin{equation}\label{estimate a}
	\begin{aligned}
		\|a_{(\xi)}\|_{C_tL^2}\leq\frac{M}{4|\Lambda|}\left(2\|\mathring{R}_q\|_{C_{[t-1,t+1]}L^1}+\frac12\delta_{q+1}\bar e\right)^{1/2},
	\end{aligned}
\end{equation}
where   $M$ denotes the universal constant from Lemma~\ref{geometric}.
Moreover,  we could get the following $C^N_{t,x}$-norm of $a_{(\xi)}$. Since the calculation is similar as in  \cite{HZZ21markov} except the explicit dependence on  $\|\mathring{R}_q\|_{C_{[t-1,t+1]}L^1}$, we put the main part in Appendix~\ref{ap:B}. In particular, by \eqref{estimate aN1}-\eqref{estimate aN01} we obtain  for $N\geq1$
\begin{equation}\label{estimate aN}
	\begin{aligned}
		\|a_{(\xi)}\|_{C^N_{t,x}}& \lesssim \ell^{-7-6N}(\|\mathring{R}_q\|_{C_{[t-1,t+1]}L^1}+1)^{N+1},
	\end{aligned}
\end{equation}
and
\begin{equation}\label{estimate aN0}
	\begin{aligned}
		\|a_{(\xi)}\|_{C^0_{t,x}}& \lesssim \ell^{-2}(\|\mathring{R}_q\|_{C_{[t-1,t+1]}L^1}+1)^{1/2}.
	\end{aligned}
\end{equation}
Here, the implicit constant depends on $N$ and in the following we only used $N\leq9$.

\medskip

With these preparations in hand,  we define the principal part $w_{q+1}^{(p)}$ of the perturbation $w_{q+1}$ as
\begin{equation}\label{principle}
	w_{q+1}^{(p)}:=\sum_{\xi\in\Lambda} a_{(\xi)}W_{(\xi)}.
\end{equation}
Since the coefficients $a_{(\xi)}$ are $(\mathcal{F}_t)_{t\geq0}$-adapted and $W_{(\xi)}$ is a deterministic function we deduce that
$w_{q+1}^{(p)}$ is also $(\mathcal{F}_t)_{t\geq0}$-adapted.
Moreover, according to (\ref{cancellation}) and (\ref{Wxi}) it follows that
\begin{equation}\label{can}w_{q+1}^{(p)}\otimes w_{q+1}^{(p)}+\mathring{R}_\ell=\sum_{\xi\in \Lambda}a_{(\xi)}^2 \mathbb{P}_{\neq0}(W_{(\xi)}\otimes W_{(\xi)})+\rho \Id,
\end{equation}
where we use the notation $\mathbb{P}_{\neq0}f:=f-\frac1{(2\pi)^3}\int_{\mT^{3}} f\dif x$.

We also define the incompressibility corrector by
\begin{equation}\label{incompressiblity}
	w_{q+1}^{(c)}:=\sum_{\xi\in \Lambda}\textrm{curl}(\nabla a_{(\xi)}\times V_{(\xi)})+\nabla a_{(\xi)}\times \textrm{curl}V_{(\xi)}+a_{(\xi)}W_{(\xi)}^{(c)},\end{equation}
with $W_{(\xi)}^{(c)}$ and $V_{(\xi)}$ being given in (\ref{corrector}).
Since $a_{(\xi)}$ is $(\mathcal{F}_t)_{t\geq0}$-adapted and $W_{(\xi)}, W_{(\xi)}^{(c)}$ and $V_{(\xi)}$ are  deterministic it follows that
$w_{q+1}^{(c)}$ is also $(\mathcal{F}_t)_{t\geq0}$-adapted.
By a direct computation we deduce that
\begin{equation*}
	w_{q+1}^{(p)}+w_{q+1}^{(c)}=\sum_{\xi\in\Lambda}\textrm{curl}\,\textrm{curl}(a_{(\xi)}V_{(\xi)}),
\end{equation*}
hence
\begin{equation*}\div(w_{q+1}^{(p)}+w_{q+1}^{(c)})=0.\end{equation*}
Next, we introduce the temporal corrector
\begin{equation}\label{temporal}w_{q+1}^{(t)}:=-\frac{1}{\mu}\sum_{\xi\in \Lambda}\mathbb{P}\mathbb{P}_{\neq0}\left(a_{(\xi)}^2\phi_{(\xi)}^2\psi_{(\xi)}^2\xi\right),\end{equation}
where $\mathbb{P}$ is the Helmholtz projection. Similarly as above, $w_{q+1}^{(t)}$ is $(\mathcal{F}_t)_{t\geq0}$-adapted and by a direct computation (see \cite[(7.20)]{BV19})  we obtain
\begin{equation}\label{equation for temporal}
	\aligned
	&\partial_t w_{q+1}^{(t)}+\sum_{\xi\in\Lambda}\mathbb{P}_{\neq0}\left(a_{(\xi)}^2\div(W_{(\xi)}\otimes W_{(\xi)})\right)
	\\
	&\qquad= -\frac{1}{\mu}\sum_{\xi\in\Lambda}\mathbb{P}\mathbb{P}_{\neq0}\partial_t\left(a_{(\xi)}^2\phi_{(\xi)}^2\psi_{(\xi)}^2\xi\right)
	+\frac{1}{\mu}\sum_{\xi\in\Lambda}\mathbb{P}_{\neq0}\left( a^2_{(\xi)}\partial_t(\phi^2_{(\xi)}\psi^2_{(\xi)}\xi)\right)
	\\&\qquad= (\Id-\mathbb{P})\frac{1}{\mu}\sum_{\xi\in\Lambda}\mathbb{P}_{\neq0}\partial_t\left(a_{(\xi)}^2\phi_{(\xi)}^2\psi_{(\xi)}^2\xi\right)
	-\frac{1}{\mu}\sum_{\xi\in\Lambda}\mathbb{P}_{\neq0}\left(\partial_t a^2_{(\xi)}(\phi^2_{(\xi)}\psi^2_{(\xi)}\xi)\right).
	\endaligned
\end{equation}
Note that the first term on the right hand side can be viewed as a pressure term $\nabla p_{1}$.

Finally, the total perturbation $w_{q+1}$ is defined by
\begin{equation}\label{wq}w_{q+1}:=w_{q+1}^{(p)}+w_{q+1}^{(c)}+w_{q+1}^{(t)},\end{equation}
which is mean zero, divergence free and $(\mathcal{F}_t)_{t\geq0}$-adapted.
The new velocity $v_{q+1}$ is defined as
\begin{equation}\label{vq}
	v_{q+1}:=v_\ell+w_{q+1}.
\end{equation}
Thus, it is also $(\mathcal{F}_t)_{t\in\mR}$-adapted.

\subsubsection{Inductive estimates for $v_{q+1}$}
\label{sss:v}
Next, we verify the inductive estimates (\ref{inductionv}) on the level $q+1$ for $v$ and we prove (\ref{iteration}).

In the following we use \cite[Lemma B.1]{CL20}.
This result is applied  to bound $w_{q+1}^{(p)}$ in $L^{2}$ whereas for the other $L^{p}$-norms we use a different approach.  By (\ref{estimate a}), \eqref{estimate a1} and (\ref{estimate aN}) we obtain
\begin{align}\label{estimate wqp}
	\|w_{q+1}^{(p)}\|_{C_tL^2}&\lesssim  \sum_{\xi\in\Lambda}\|a_{(\xi)}\|_{{C_{t}}L^2}\|W_{(\xi)}\|_{{C_{t}}L^2}+\frac1{(\lambda_{q+1} r_\bot)^{1/2}}\|a_{(\xi)}\|_{C^1_{{t,x}}}\|W_{(\xi)}\|_{{C_{t}}L^2}\nonumber
	\\&\leq\frac{M_0}8(\|\mathring{R}_q\|_{C_{[t-1,t+1]}L^1}+\bar e\delta_{q+1})^{1/2}+\frac1{(\lambda_{q+1} r_\perp)^{1/2}}\ell^{-13}(\|\mathring{R}_q\|_{C_{[t-1,t+1]}L^1}+1)^{2}
\end{align}
where we used  the fact that due to (\ref{intermittent}) together with the normalizations \eqref{eq:phi}, \eqref{eq:psi} it holds  $\|W_{(\xi)}\|_{L^2}\simeq 1$ uniformly in all the involved parameters. Here, we may  choose $M_0=cM\geq1$ with a universal constant $c$.

For a general $L^p$-norm we apply (\ref{bounds}) and (\ref{estimate aN0}) to deduce for  $p\in(1,\infty)$
\begin{equation}\label{principle est1}
	\aligned
	\|w_{q+1}^{(p)}\|_{C_tL^p}&\lesssim \sum_{\xi\in \Lambda}\|a_{(\xi)}\|_{C^0_{t,x}}\|W_{(\xi)}\|_{C_tL^p}\lesssim \ell^{-2}(\|\mathring{R}_q\|_{C_{[t-1,t+1]}L^1}+1)^{1/2}r_\perp^{2/p-1}r_\|^{1/p-1/2},
	\endaligned
\end{equation}
\begin{equation}\label{correction est}
	\aligned
	\|w_{q+1}^{(c)}\|_{C_tL^p}&\lesssim\sum_{\xi\in \Lambda}\left(\|a_{(\xi)}\|_{C^0_{t,x}}\|W_{(\xi)}^{(c)}\|_{C_tL^p}+\|a_{(\xi)}\|_{C^2_{t,x}}\|V_{(\xi)}\|_{C_tW^{1,p}}\right)
	\\&\lesssim \ell^{-19}(\|\mathring{R}_q\|_{C_{[t-1,t+1]}L^1}+1)^3r_\perp^{2/p-1}r_\|^{1/p-1/2}\left(r_\perp r_\|^{-1}+\lambda_{q+1}^{-1}\right)
	\\&\lesssim \ell^{-19}(\|\mathring{R}_q\|_{C_{[t-1,t+1]}L^1}+1)^3r_\perp^{2/p}r_\|^{1/p-3/2},
	\endaligned
\end{equation}
and
\begin{equation}\label{temporal est1}
	\aligned
	\|w_{q+1}^{(t)}\|_{C_tL^p}&\lesssim \mu^{-1}\sum_{\xi\in\Lambda}\|a_{(\xi)}\|_{C^0_{t,x}}^2\|\phi_{(\xi)}\|_{L^{2p}}^2\|\psi_{(\xi)}\|_{C_tL^{2p}}^2
	\\
	&\lesssim\ell^{-4}(\|\mathring{R}_q\|_{C_{[t-1,t+1]}L^1}+1)r_\perp^{2/p-1}r_\|^{1/p-2}(\mu^{-1}r_\perp^{-1}r_\|)
	\\&= \ell^{-4}(\|\mathring{R}_q\|_{C_{[t-1,t+1]}L^1}+1)r_\perp^{2/p-1}r_\|^{1/p-2}\lambda_{q+1}^{-1}.
	\endaligned
\end{equation}
We note that for $p=2$ \eqref{principle est1} provides a worse bound than \eqref{estimate wqp}.
Hence, we obtain for $p=\frac{32}{32-7\alpha}>1$ so that  $r_\perp^{2/p-2}r_\|^{1/p-1}\leq \lambda_{q+1}^\alpha$ it holds that
\begin{equation}\label{corr temporal}
	\aligned
	&\|w_{q+1}\|_{C_tL^p}\lesssim  \lambda_{q+1}^{-8/7+5\alpha}(\|\mathring{R}_q\|_{C_{[t-1,t+1]}L^1}+1)^3,
	\endaligned
\end{equation}
where we use (\ref{ell}) and the fact that $\lambda_{q+1}^{19\alpha-\frac{1}{7}}<1$ by our choice of $\alpha$. The bound \eqref{corr temporal} will be used below in the estimation of the Reynolds stress.

Combining (\ref{estimate wqp}), (\ref{correction est}) and (\ref{temporal est1}) we obtain
\begin{equation}\label{estimate wq}
	\aligned
	\|w_{q+1}\|_{C_tL^2}&\leq \frac{M_0}{8}\|\mathring{R}_q\|^{1/2}_{C_{[t-1,t+1]}L^1}+\bar e^{1/2}\delta_{q+1}^{1/2}\frac{M_0}{4}+\lambda_{q+1}^{-1/{14}}\ell^{-13}(\|\mathring{R}_q\|_{C_{[t-1,t+1]}L^1}+1)^{2}
	\\&\quad+\lambda_{q+1}^{-2/7}\ell^{-19}(\|\mathring{R}_q\|_{C_{[t-1,t+1]}L^1}+1)^3+\lambda_{q+1}^{-1/7}\ell^{-4}(\|\mathring{R}_q\|_{C_{[t-1,t+1]}L^1}+1)
	\\&\leq \frac{M_0}{8}\|\mathring{R}_q\|^{1/2}_{C_{[t-1,t+1]}L^1}+\bar e^{1/2}\delta_{q+1}^{1/2}\frac{M_0}{4}+\lambda_{q+1}^{-1/{14}+26\alpha}(\|\mathring{R}_q\|_{C_{[t-1,t+1]}L^1}^3+1)\frac{M_0}8.
	\endaligned
\end{equation}
Thus by \eqref{eq:R} and \eqref{bd:R} we obtain
\begin{equation}\label{estimate wq1}
	\aligned
	\$w_{q+1}\$_{L^2,2r} &\leq\frac{M_0}{4}\$\mathring{R}_q\$^{1/2}_{L^1,r}+\bar e^{1/2}\delta_{q+1}^{1/2}\frac{M_0}{4}+\frac{M_0}4\lambda_{q+1}^{-1/{14}+26\alpha}(\$\mathring{R}_q\$^3_{L^1,6r}+1)
	\\&\leq\frac{M_0}{4}\bar e^{1/2}\delta_{q+1}^{1/2}+\frac14M_0\bar e^{1/2}\delta_{q+1}^{1/2}+\frac{M_0}{4}\lambda_{q+1}^{-1/{14}+26\alpha}{\lambda^6_q}((6^qr\cdot 24L^2)^{3 (6^q)}+1)
	\\&\leq\frac{M_0}{2}\bar e^{1/2}\delta_{q+1}^{1/2}+\frac{M_0}4\lambda_{q+1}^{-1/{14}+26\alpha}(\lambda_q^{{9}}+1)
	\\&\leq\frac{3M_0}4\bar e^{1/2}\delta_{q+1}^{1/2}.
	\endaligned
\end{equation}
Here we used $a>(144 L^2r)^c$ and $-1/{14}+26\alpha+{9}/b<-\beta$ by the choice of parameters in Section~\ref{s:par}.
The bound \eqref{estimate wq1} can be directly combined with \eqref{error ps} and the definition of the velocity $v_{q+1}$ \eqref{vq} to deduce
$$
\$v_{q+1}-v_{q}\$_{L^{2},2r}\leq \$w_{q+1}\$_{L^{2},2r}+\$v_{\ell}-v_{q}\$_{L^{2},2r}\leq M_0\bar e^{1/2}\delta_{q+1}^{1/2},
$$
hence \eqref{iteration} holds and \eqref{inductionv} follows at the level of $q+1$.
Moreover,  by \eqref{estimate wq} and \eqref{bd:R} we obtain
\begin{equation}\label{estimate wq2}
	\aligned
	\$w_{q+1}\$_{L^2,m} &\leq\frac{M_0}{4}\$\mathring{R}_q\$^{1/2}_{L^1,m/2}+\bar e^{1/2}\delta_{q+1}^{1/2}\frac{M_0}{4}+\frac{M_0}4\lambda_{q+1}^{-1/{14}+26\alpha}(\$\mathring{R}_q\$^3_{L^1,3m}+1)
	\\&\leq\frac{M_0}{4}(6^q2mL^2)^{\frac126^q}{\lambda_q}+M_0\bar e^{1/2}\delta_{q+1}^{1/2}+\frac{M_0}4\lambda_{q+1}^{-1/{14}+26\alpha}(6^q\cdot 12mL^2)^{3 (6^q)}{\lambda^6_q}
	\\&\leq\frac{M_0}2(6^q\cdot 12mL^2)^{3 (6^q)}{\lambda_q}+M_0\bar e^{1/2}\delta_{q+1}^{1/2}.
	\endaligned
\end{equation}
Thus combined with \eqref{eq:vl} and \eqref{inductionv m}  we deduce that \eqref{inductionv m} holds at level of $q+1$.

As the next step, we shall verify the first bound in \eqref{inductionv C1}. The $C^1_{t,x}$-bound follows similarly as in \cite{HZZ21markov} but with an explicit dependence on $\|\mathring{R}_q\|_{C_{[t-1,t+1]}L^1}$. Hence, we omit most details for these estimates.
Using (\ref{estimate aN})-\eqref{estimate aN0} and (\ref{bounds}) we have
\begin{equation}\label{principle est2}
	\aligned
	\|w_{q+1}^{(p)}\|_{C^1_{t,x}}
	&\lesssim \ell^{-13}(\|\mathring{R}_q\|_{C_{[t-1,t+1]}L^1}+1)^{2}r_\perp^{-1}r_\|^{-1/2}\lambda_{q+1}^2,
	\endaligned
\end{equation}
\begin{equation}\label{correction est2}
	\aligned
	\|w_{q+1}^{(c)}\|_{C^1_{t,x}}
	&\lesssim\ell^{-25}(\|\mathring{R}_q\|_{C_{[t-1,t+1]}L^1}+1)^4r_\|^{-3/2}\lambda_{q+1}^2,
	\endaligned
\end{equation}
and
\begin{equation}\label{temporal est2}
	\aligned
	\|w_{q+1}^{(t)}\|_{C^1_{t,x}}&\leq \frac{1}{\mu}\sum_{\xi\in\Lambda}[\|a^2_{(\xi)}
	\phi^2_{(\xi)}\psi^2_{(\xi)}\|_{C_tW^{1+\alpha,p}}+\|a^2_{(\xi)}\phi^2_{(\xi)}\psi^2_{(\xi)}\|_{C^1_tW^{\alpha,p}}]
	\\&\lesssim\ell^{-21}(\|\mathring{R}_q\|_{C_{[t-1,t+1]}L^1}+1)^{3}r_\perp^{-1}r_\|^{-2}\lambda_{q+1}^{1+\alpha},
	\endaligned
\end{equation}
where we chose $p$ large enough and applied the Sobolev embedding in the first inequality. This was needed because $\mathbb{P}\mathbb{P}_{\neq0}$ is not a bounded operator on $C^0$. In the last inequality in \eqref{temporal est2}, we used interpolation and an extra $\lambda_{q+1}^\alpha$ appeared.
Combining   \eqref{principle est2}, \eqref{correction est2}, \eqref{temporal est2} with (\ref{ell}) we obtain
\begin{equation*}
	\aligned
	\|v_{q+1}\|_{C^1_{t,x}}&\leq \|v_\ell\|_{C^1_{t,x}}+\|w_{q+1}\|_{C^1_{t,x}}\\
	&\leq (\|\mathring{R}_q\|_{C_{[t-1,t+1]}L^1}+1)^{4}\left(C\lambda_{q+1}^{26\alpha+22/7}+C\lambda_{q+1}^{50\alpha+20/7}+C\lambda_{q+1}^{43\alpha+3}\right)+\|v_q\|_{C^1_{[t-1,t+1],x}}.
	\endaligned
\end{equation*}
Thus,
\begin{equation*}
	\aligned
	\$v_{q+1}\$_{C^1_{t,x},m}&\lesssim (\$\mathring{R}_q\$^4_{L^1,4m}+1)\lambda_{q+1}^{26\alpha+22/7}+\$v_q\$_{C^1_{t,x},m}
	\\&\leq \lambda_{q+1}^{23/7}(6^q\cdot 16mL^2)^{4(6^q)},
	\endaligned
\end{equation*}
which implies the first inequality in \eqref{inductionv C1} holds true on the level $q+1$.

Similarly, we have
\begin{equation}\label{principle est2-c2}
	\aligned
	\|w_{q+1}^{(p)}\|_{C^2_{t,x}}
	&\lesssim \ell^{-19}(\|\mathring{R}_q\|_{C_{[t-1,t+1]}L^1}+1)^3r_\perp^{-1}r_\|^{-1/2}\lambda_{q+1}^2\left(1+\frac{r_\perp \mu}{r_\|}\right)^2
	\\&\lesssim \lambda_{q+1}^{38\alpha+{36}/7}(\|\mathring{R}_q\|_{C_{[t-1,t+1]}L^1}+1)^3,
	\endaligned
\end{equation}
\begin{equation}\label{correction est2-c2}
	\aligned
	\|w_{q+1}^{(c)}\|_{C^2_{t,x}}
	&\lesssim \ell^{-31}(\|\mathring{R}_q\|_{C_{[t-1,t+1]}L^1}+1)^5r_\|^{-3/2}\Big(\lambda_{q+1}\frac{r_\perp \mu}{r_\|}\Big)^2
	\\&\lesssim \lambda_{q+1}^{{34}/7+62\alpha}(\|\mathring{R}_q\|_{C_{[t-1,t+1]}L^1}+1)^5,
	\endaligned
\end{equation}
and
\begin{equation}\label{temporal est2-c2}
	\aligned
	\|w_{q+1}^{(t)}\|_{C^2_{t,x}}&\lesssim \ell^{-27}(\|\mathring{R}_q\|_{C_{[t-1,t+1]}L^1}+1)^{5}r_\perp^{-2}r_\|^{-1}\lambda_{q+1}^{2+\alpha}\mu^{-1}\left(1+\frac{r_\perp \mu}{r_\|}\right)^2
	\\&\lesssim(\|\mathring{R}_q\|_{C_{[t-1,t+1]}L^1}+1)^{5}\lambda_{q+1}^{55\alpha+5}.
	\endaligned
\end{equation}
Hence, we obtain
\begin{align*}
	\|v_{q+1}\|_{C^2_{t,x}}\lesssim (\|\mathring{R}_q\|_{C_{[t-1,t+1]}L^1}+1)^{5}\lambda_{q+1}^{{38\alpha+{36}/7}}+\|v_{q}\|_{C^2_{[t-1,t+1],x}},
\end{align*}
which implies
\begin{align*}
	\$v_{q+1}\$_{C^2_{t,x},m}&\lesssim (\$\mathring{R}_q\$_{C_{[t-1,t+1]}L^1,5m}^5+1)\lambda_{q+1}^{{38\alpha+{36}/7}}+\$v_{q}\$_{C^2_{t,x},m}
	\\&\leq\lambda_{q+1}^{{37}/7}\cdot (6^q\cdot 20 mL^2)^{5\cdot 6^{q}}.
\end{align*}
Hence, the second inequality in \eqref{inductionv C1} holds true on the level $q+1$.

We conclude this part with further estimates of the perturbations $w^{(p)}_{q+1}$, $w^{(c)}_{q+1}$ and $w^{(t)}_{q+1}$, which will be used below in order to bound the Reynolds stress $\mathring{R}_{q+1}$ and to establish \eqref{bdvq:w1p}.
These estimates follow similarly as in \cite{HZZ21markov} with an explicit dependence on $\|\mathring{R}_q\|_{C_{[t-1,t+1]}L^1}$.  We omit most details and derive  the following estimates: by using (\ref{ell}), (\ref{estimate aN}), \eqref{estimate aN0} and (\ref{bounds})
\begin{equation}\label{principle est22}
	\aligned
	&\|w_{q+1}^{(p)}+w_{q+1}^{(c)}\|_{C_tW^{1,p}}
	\\&\leq\sum_{\xi\in\Lambda}
	\|\textrm{curl\,}\textrm{curl}(a_{(\xi)}V_{(\xi)})\|_{C_tW^{1,p}}
	\\
	&\lesssim r_\perp^{2/p-1}r_\|^{1/p-1/2}\left(\ell^{-25}
	\lambda_{q+1}^{-2}+\ell^{-19}\lambda_{q+1}^{-1}+\ell^{-13}+\ell^{-2}\lambda_{q+1}\right)(\|\mathring{R}_q\|_{C_{[t-1,t+1]}L^1}+1)^{4}\\
	&\lesssim r_\perp^{2/p-1}r_\|^{1/p-1/2}\ell^{-2}\lambda_{q+1}(\|\mathring{R}_q\|_{C_{[t-1,t+1]}L^1}+1)^{4},
	\endaligned
\end{equation}
and
\begin{equation}\label{corrector est2}
	\aligned
	\|w_{q+1}^{(t)}\|_{C_tW^{1,p}}
	&\lesssim \frac{1}{\mu}r_\perp^{2/p-2}r_\|^{1/p-1}\left(\ell^{-{15}}+\ell^{-4}\lambda_{q+1}\right)(\|\mathring{R}_q\|_{C_{[t-1,t+1]}L^1}+1)^{5/2}
	\\&\lesssim r_\perp^{2/p-2}r_\|^{1/p-1}\ell^{-4}\lambda_{q+1}^{-2/7}(\|\mathring{R}_q\|_{C_{[t-1,t+1]}L^1}+1)^{{5/2}}.
	\endaligned
\end{equation}

We then obtain for $p=\frac{32}{32-7\alpha}$
\begin{align}
	\|w_{q+1}\|_{C_tW^{1,p}}&\lesssim (r_\perp^{2/p-1}r_{\|}^{1/p-1/2}\ell^{-2}\lambda_{q+1}+r_\perp^{2/p-2}r_{\|}^{1/p-1}\ell^{-4}\lambda_{q+1}^{-2/7})(1+\|R_q\|_{C_{[t-1,t+1]}L^1})^4\no\\
	&\lesssim (\lambda_{q+1}^{-2/7+9\alpha}+\lambda_{q+1}^{-1/7+5\alpha})(1+\|R_q\|_{C_{[t-1,t+1]}L^1})^4.\label{eq:w1p}
\end{align}
Taking expectation we obtain
\begin{align*}
	\E\|w_{q+1}\|_{C_tW^{1,p}}^r\lesssim \lambda_{q+1}^{-{r}/7+5r\alpha}(1+\E\|R_q\|_{C_tL^1}^{4r})\lesssim \lambda_{q+1}^{-{r}/7+5r\alpha}(6^q\cdot 16r L^2)^{6^q\cdot 4r}{\lambda^{8r}_q}\leq \lambda_{q+1}^{-{r}/7+{8}r\alpha},
\end{align*}
where we used $(6^q\cdot 16 rL^2)^{6^q}\leq \lambda_q$ , $\lambda_q^4\leq \lambda_{q+1}^\alpha$ and we chose $a$ large enough to absorb the constant.
Moreover, by \eqref{vqc1} we obtain
\begin{align*}
	\E\|v_\ell-v_{q}\|_{C_{t}W^{1,p}}^r\leq \ell ^r \lambda_q^{6r}\leq {\frac12}\lambda_{q+1}^{-\alpha r/2}.
\end{align*}
Now, \eqref{bdvq:w1p} follows.

\subsubsection{Proof of  \eqref{p:gamma}}\label{s:en}

We define
$$\delta E(t):=\Big| e(t)(1-\delta_{q+2})-\E\|v_{q+1}(t)+{z_{q+1}}(t)\|_{L^2}^2\Big|.$$

\bp It holds for $t\in\mR$
\begin{align}
	\delta E(t)\leq& \frac14\delta_{q+2}e(t).
\end{align}
\ep
\begin{proof}
	By definition of $\gamma_q$ we find
	\begin{align}\label{eq:deltaE}
		\begin{aligned}
			\delta E(t)&\leq \E\big|\|w_{q+1}^{(p)}\|_{L^2}^2-3\gamma_q(2\pi)^3\big|+\E\|w_{q+1}^{(c)}+w_{q+1}^{(t)}\|_{L^2}^2+2\E\|(v_\ell+z_{q+1})\cdot(w_{q+1}^{(c)}+w_{q+1}^{(t)})\|_{L^1}\\
			&\qquad+2\E\|(v_\ell+z_{q+1})\cdot w_{q+1}^{(p)}\|_{L^1}+2\E\|w_{q+1}^{(p)}\cdot(w_{q+1}^{(c)}+w_{q+1}^{(t)})\|_{L^1}\\
			&\qquad+\E\|v_\ell-v_q+z_{q+1}-z_q\|_{L^2}^2+2\E\|(v_\ell-v_q+z_{q+1}-z_q)\cdot(v_q+z_q)\|_{L^1},
		\end{aligned}
	\end{align}
	which shall be estimated.
	Let us begin with the bound of the first term on the right hand side of \eqref{eq:deltaE}. We use \eqref{can} and the fact that $\mathring{R}_\ell$ is traceless to deduce for $t\in\mR$
	\begin{align*}
		|  w_{q+1}^{(p)}|^2-3\gamma_{q}=6\sqrt{\ell^2+|\mathring{R}_\ell|^2}+3(\gamma_\ell-\gamma_q)+\sum_{\xi\in \Lambda}a_{(\xi)}^2\mathbb{P}_{\neq0}|W_{(\xi)}|^2
		,
	\end{align*}
	hence
	\begin{align}\label{eq:gg ps}
	\begin{aligned}
&		\mathbf{E}|\|  w_{q+1}^{(p)}\|_{L^2}^2-3\gamma_{q+1}(2\pi)^3|\\
		&\qquad\leq 6\cdot(2\pi)^3\ell+6\mathbf{E}\|\mathring{R}_\ell\|_{L^1}+3\cdot(2\pi)^3|\gamma_\ell-\gamma_q|+\mathbf{E}\sum_{\xi\in \Lambda}\Big|\int a_{(\xi)}^2\mathbb{P}_{\neq0}|W_{(\xi)}|^2\Big|.
		\end{aligned}
	\end{align}
	Here we estimate each term separately. Using \eqref{ell1}  we find
	\begin{align*}
		6\cdot(2\pi)^3\ell\leq 6\cdot(2\pi)^3\lambda_{q+1}^{-{3\alpha}/2}\leq  \frac{1}{ 48}\lambda_{q+1}^{-2\beta b}e(t)\leq\frac{1}{ 48}\delta_{q+2}e(t),
	\end{align*}
	which requires $2\beta b<{3\alpha}/2$ and choosing $a$ large to absorb the constant.
	Using \eqref{eq:R} on $\mathring{R}_q$ and $\supp\varphi_\ell\subset [0,\ell]$ we know for $t\in\mR$
	\begin{align*}
		6\mathbf{E}\|\mathring{R}_\ell(t)\|_{L^1}\leq  \frac18\delta_{q+2}e(t).
	\end{align*}
	For the third term in \eqref{eq:gg ps} we use \eqref{eq:u0} and \eqref{inductionv},   \eqref{vqc1} to have for $0\leq \delta\leq 1/6$
	\begin{align*}
		3\cdot(2\pi)^3|\gamma_\ell-\gamma_q|&\lesssim \ell \|e'\|_{C^{0}_{t-1}}+\ell  \E\|v_q\|_{C^1_{t-1,x}}(\|v_q\|_{C_{t-1}L^2}+{\|z_q\|_{C_{t-1}L^2}})\\
		&\qquad+\ell^{1/2-\delta}\E\|{z_q}\|_{C_{t-1}^{1/2-\delta}L^2}(\|v_q\|_{C_{t-1}L^2}+\|{z_q}\|_{C_{t-1}L^2})
		\\
		&\lesssim \ell \tilde e+\ell \lambda_q^{4}(\bar{e}^{1/2}M_0+L)+\ell^{1/2-\delta}L(M_0\bar{e}^{1/2}+L)\\
		&\lesssim \lambda_{q+1}^{-3\alpha/2}\tilde{e}+\lambda_{q+1}^{-\alpha} (M_0\bar{e}^{1/2}+L)+\lambda_{q+1}^{-\frac{3\alpha}{2}(1/2-\delta)}L (M_0\bar{e}^{1/2}+L)
		\\&\lesssim  \lambda_{q+1}^{-{\alpha}/{2}} (\bar{e}+\tilde e+L^2)\leq \frac{1}{ 48}\delta_{q+2}e(t),
	\end{align*}
	where we choose $a$ large to absorb the constant.

	For the last term in \eqref{eq:gg ps}
	we apply \eqref{estimate aN}, \eqref{estimate aN0} and $\|a_{(\xi)}^2\|_{C^N}\lesssim \|a_{(\xi)}\|_{C^0}\|a_{(\xi)}\|_{C^N}$ to bound
	\begin{align*}
		\sum_{\xi\in\Lambda}&\Big|\int_{\mathbb{T}^{3}} a_{(\xi)}^2\mathbb{P}_{\neq0}|W_{(\xi)}|^2 \dif x\Big|=\sum_{\xi\in\Lambda}\Big|\int_{\mathbb{T}^{3}} a_{(\xi)}^2\mathbb{P}_{\geq r_\perp\lambda_{q+1}/2}|W_{(\xi)}|^2\dif x\Big|\\
		&=\sum_{\xi\in\Lambda}\Big|\int_{\mathbb{T}^{3}} |\nabla|^Na_{(\xi)}^2|\nabla|^{-N}\mathbb{P}_{\geq r_\perp\lambda_{q+1}/2}|W_{(\xi)}|^2 \dif x\Big|
		\\
		&\lesssim\sum_{\xi\in\Lambda}\|a_{(\xi)}^2\|_{C^N}(r_\perp\lambda_{q+1})^{-N}\||W_{(\xi)}|^2\|_{L^2}\lesssim\ell^{-6N-9}(\|\mathring{R}_q\|_{C_{[t-1,t+1]}L^1}+1)^{N+3/2}(r_\perp\lambda_{q+1})^{-N}r_\perp^{-1}r_{\|}^{-\frac12}\\
		&\leq\ell^{-6N-9}(\|\mathring{R}_q\|_{C_{[t-1,t+1]}L^1}+1)^{N+3/2}\lambda_{q+1}^{\frac{8-N}7}.
	\end{align*}
	Thus
	\begin{align*}
		\mathbf{E}\sum_{\xi\in\Lambda}\Big|\int_{\mathbb{T}^{3}} a_{(\xi)}^2\mathbb{P}_{\neq0}|W_{(\xi)}|^2 \dif x\Big|
		&\lesssim \lambda_{q+1}^{(12N+18)\alpha+\frac{8-N}7 }(6^q\cdot4(N+3/2)L^2)^{6^q(N+3/2)}{\lambda_q^{2N+3}}
		\\&\leq\lambda_{q+1}^{{132}\alpha-1/7}\leq\frac1{48}\delta_{q+2}e(t).
	\end{align*}
	Here we may choose $N=9$, $a>[252L^2]^{3c}$ such that $(6^q\cdot4(N+3/2)L^2)^{6^q(N+3/2)}<\lambda_q^4<\lambda_{q+1}^\alpha$ and   use $2\beta b <1/{7}-{132}\alpha$, $\lambda_q^{21}<\lambda_{q+1}^{5\alpha}$.
	This completes the bound for \eqref{eq:gg ps}.

	Going back to \eqref{eq:deltaE}, it remains to control
	\begin{align*}
		&\E\|w_{q+1}^{(c)}+w_{q+1}^{(t)}\|_{L^2}^2+2\E\|(v_\ell+z_{q+1})(w_{q+1}^{(c)}+w_{q+1}^{(t)})\|_{L^1}\\
		&\qquad+2\E\|(v_\ell+z_{q+1})w_{q+1}^{(p)}\|_{L^1}+2\E\|w_{q+1}^{(p)}(w_{q+1}^{(c)}+w_{q+1}^{(t)})\|_{L^1}\\
		&\qquad+\E\|v_\ell-v_q+z_{q+1}-z_q\|_{L^2}^2+2\E\|(v_\ell-v_q+z_{q+1}-z_q)(v_q+z_q)\|_{L^1}.
	\end{align*}
	Using  the estimates \eqref{correction est}, \eqref{temporal est1} and \eqref{ell} we  have
	\begin{align*}
		\mathbf{E}\|  w_{q+1}^{(c)}+  w_{q+1}^{(t)}\|_{L^2}^2&
		\lesssim (1+\$\mathring{R}_q\$_{L^1,6}^6)\lambda_{q+1}^{{76}\alpha-4/7}
		+(1+\$\mathring{R}_q\$_{L^1,2}^2)\lambda_{q+1}^{16\alpha-2/7}
		\\&\lesssim(6^q\cdot24L^2)^{6^{q+1}}{\lambda^{12}_q} \lambda_{q+1}^{16\alpha-2/7}\leq \frac{\delta_{q+2}}{48}e(t),
	\end{align*}
	where we use a  similar bound for the parameters as above. Next, we use \eqref{eq:vl} together with \eqref{estimate wqp} to have
	\begin{align*}
		&\mathbf{E}[2\|(v_\ell+z_{q+1})(  w_{q+1}^{(c)}+  w_{q+1}^{(t)})\|_{L^1}+2\|  w_{q+1}^{(p)}(  w_{q+1}^{(c)}+  w_{q+1}^{(t)})\|_{L^1}]
		\lesssim   {(M_{0}\bar{e}^{1/2}+L)}\$  w_{q+1}^{(c)}+  w_{q+1}^{(t)}\$_{L^2,2}\\
		&\qquad \lesssim {(M_{0}\bar{e}^{1/2}+L)} (1+\$\mathring{R}_q\$_{L^1,6}^3)\lambda_{q+1}^{8\alpha-1/7}
		\\
		&\qquad\leq \frac{1}{48}
		\lambda_{q+1}^{-2\beta b}\leq \frac{\delta_{q+2}}{48}e(t),
	\end{align*}
	with similar arguments for the second  last inequality as above.
	We employ \eqref{ell},  \eqref{vqc1}, \eqref{principle est1} as well as $\|v_\ell\|_{C^1_{t,x}}\leq \|v_q\|_{C^1_{[t-1,t+1],x}} $ to have for every $\eps>0$
	\begin{align}\label{eq:00}
		\begin{aligned}
			2\mathbf{E}\|(v_\ell +z_{q+1})  w_{q+1}^{(p)}\|_{L^1}
			&\lesssim (\$v_\ell\$_{L^\infty,2}+\$z_q\$_{L^\infty,2})\$  w_{q+1}^{(p)}\$_{L^1,2}+\$z_{q+1}-z_q\$_{L^4,2}\$  w_{q+1}^{(p)}\$_{L^{4/3},2}
			\\&\lesssim (\lambda_q^{4}+\lambda_{q+1}^{\alpha/8}L)\ell^{-2}\delta_{q+2}^{1/2}r_{\perp}^{1-\eps}r_{\|}^{\frac12(1-\eps)}+\lambda_{q+1}^{-\frac\alpha8(\frac14-\delta)}L\ell^{-2}\delta_{q+2}^{1/2}r_{\perp}^{1/2}r_{\|}^{1/4}\\
			&\lesssim L\lambda_{q+1}^{5\alpha-\frac87(1-\eps)}+L\lambda_{q+1}^{4\alpha-\frac47}\leq \frac{1}{96}
			\lambda_{q+1}^{-2\beta b}\leq \frac{\delta_{q+2}}{96}e(t).
		\end{aligned}
	\end{align}
	
	For the last  terms, we apply \eqref{error ps} and obtain for $0<\delta<1/9$
	\begin{align}\label{eq:01}
		\begin{aligned}
			\mathbf{E}&\|v_\ell-v_q+z_{q+1}-z_q\|_{L^2}^2+2\E\|(v_\ell-v_q+z_{q+1}-z_q)(v_q+z_q)\|_{L^1}
			\\&\lesssim \$v_\ell-v_q\$_{L^2,2}(\$v_q\$_{L^2,2}+\$z_q\$_{L^2,2}+1)+\E\|z_{q+1}-z_q\|_{L^2}^2
			\\&\qquad+\$z_{q+1}-z_q\$_{L^2,2}(\$v_q\$_{L^2,2}+\$z_q\$_{L^2,2})\\
			&\lesssim \ell\lambda_q^{4}  {(M_{0} \bar e^{1/2}+L)}+\lambda_{q+1}^{-\frac\alpha8(1-\delta)}{(M_{0} \bar e^{1/2}+L)}\\&\leq \lambda_{q+1}^{-\alpha/2} ({M_{0} \bar e^{1/2}+L)}+\lambda_{q+1}^{-\frac\alpha8(1-\delta)} ({M_{0} \bar e^{1/2}+L)}\\&\leq \frac{1}{96}
			\lambda_{q+1}^{-2\beta b}\leq \frac{\delta_{q+2}}{96}e(t).
		\end{aligned}
	\end{align}
Here, we choose again $a$ large enough to absorb the extra constant.
	
	Combining the above estimates, \eqref{p:gamma} follows on the level  $q+1$.
\end{proof}

\subsubsection{Definition of the Reynolds stress $\mathring{R}_{q+1}$}\label{s:def}

Subtracting from (\ref{induction ps}) at level $q+1$ the system (\ref{mollification}), we obtain
\begin{equation}\label{stress}
	\aligned
	&\div\mathring{R}_{q+1}-\nabla p_{q+1}\\
	&=\underbrace{-\Delta w_{q+1}+\partial_t(w_{q+1}^{(p)}+w_{q+1}^{(c)})+\div((v_\ell+z_\ell)\otimes w_{q+1}+w_{q+1}\otimes (v_\ell+z_\ell))}_{\div(R_{\textrm{lin}})+\nabla p_{\textrm{lin}}}
	\\&\quad+\underbrace{\div\left((w_{q+1}^{(c)}+w_{q+1}^{(t)})\otimes w_{q+1}+w_{q+1}^{(p)}\otimes (w_{q+1}^{(c)}+w_{q+1}^{(t)})\right)}_{\div(R_{\textrm{cor}})+\nabla p_{\textrm{cor}}}
	\\&\quad+\underbrace{\div(w_{q+1}^{(p)}\otimes w_{q+1}^{(p)}+\mathring{R}_\ell)+\partial_tw_{q+1}^{(t)}}_{\div(R_{\textrm{osc}})+\nabla p_{\textrm{osc}}}
	\\&\quad+\underbrace{\div\left(v_{q+1}{\otimes}{z_{q+1}}-v_{q+1}{\otimes}z_\ell+{z_{q+1}}{\otimes}v_{q+1}-z_\ell{\otimes}v_{q+1}
		+{z_{q+1}}{\otimes}{z_{q+1}}-z_\ell{\otimes}z_\ell\right)+(z_\ell-z_{q+1})}_{\div(R_{\textrm{com}1})+\nabla p_{\textrm{com}1}}
	\\&\quad+\div(R_{\textrm{com}})-\nabla p_\ell.
	\endaligned
\end{equation}
By using $\mathcal{R}$ introduced in Section \ref{s:not} we define
\begin{equation*}
	R_{\textrm{lin}}:=-\mathcal{R}\Delta w_{q+1}+\mathcal{R}\partial_t(w_{q+1}^{(p)}+w_{q+1}^{(c)})
	+(v_\ell+z_\ell)\mathring\otimes w_{q+1}+w_{q+1}\mathring\otimes (v_\ell+z_\ell),
\end{equation*}
\begin{equation*}
	R_{\textrm{cor}}:=(w_{q+1}^{(c)}+w_{q+1}^{(t)})\mathring{\otimes} w_{q+1}+w_{q+1}^{(p)}\mathring{\otimes} (w_{q+1}^{(c)}+w_{q+1}^{(t)}),
\end{equation*}
\begin{equation*}
	R_{\textrm{com}1}:=v_{q+1}\mathring{\otimes}{z_{q+1}}-v_{q+1}\mathring{\otimes}z_\ell+{z_{q+1}}\mathring{\otimes}v_{q+1}-z_\ell\mathring{\otimes}v_{q+1}
	+{z_{q+1}}\mathring{\otimes}{z_{q+1}}-z_\ell\mathring{\otimes}z_\ell+\mathcal{R}(z_\ell-z_{q+1}).
\end{equation*}

In order to define the remaining oscillation error from the third line in \eqref{stress}, we apply (\ref{can}) and (\ref{equation for temporal}) to obtain
\begin{align*}
	&\div(w_{q+1}^{(p)}\otimes w_{q+1}^{(p)}+\mathring{R}_\ell)+\partial_tw_{q+1}^{(t)}\\
	&\quad=\sum_{\xi\in\Lambda}\div\left(a^2_{(\xi)}\mathbb{P}_{\neq0}(W_{(\xi)}\otimes W_{(\xi)})\right)+\nabla \rho+\partial_tw_{q+1}^{(t)}
	\\
	&\quad=\sum_{\xi\in\Lambda}\mathbb{P}_{\neq0}\left(\nabla a^2_{(\xi)}\mathbb{P}_{\neq0}(W_{(\xi)}\otimes W_{(\xi)})\right)+\nabla \rho+\sum_{\xi\in\Lambda}\mathbb{P}_{\neq0}\left(a^2_{(\xi)}\div(W_{(\xi)}\otimes W_{(\xi)})\right)+\partial_tw_{q+1}^{(t)}
	\\
	&\quad=\sum_{\xi\in\Lambda}\mathbb{P}_{\neq0}
	\left(\nabla a_{(\xi)}^2\mathbb{P}_{\neq0}(W_{(\xi)}\otimes W_{(\xi)}) \right)+\nabla \rho+\nabla p_1-\frac{1}{\mu}\sum_{\xi\in\Lambda}\mathbb{P}_{\neq0}
	\left(\partial_t a_{(\xi)}^2(\phi_{(\xi)}^2\psi_{(\xi)}^2\xi) \right)
\end{align*}
Therefore,
\begin{equation*}
	\aligned
	R_{\textrm{osc}}:=&\sum_{\xi\in\Lambda}\mathcal{B}
	\left(\nabla a_{(\xi)}^2,\mathbb{P}_{\neq0}(W_{(\xi)}\otimes W_{(\xi)}) \right)-\frac{1}{\mu}\sum_{\xi\in\Lambda}\mathcal{R}
	\left(\partial_t a_{(\xi)}^2(\phi_{(\xi)}^2\psi_{(\xi)}^2\xi) \right)=:R_{\textrm{osc}}^{(x)}+R_{\textrm{osc}}^{(t)},
	\endaligned
\end{equation*}
with $\cB$ given in Section \ref{s:not}.

Finally we define the Reynolds stress on the level $q+1$ by
\begin{equation*}\aligned
	\mathring{R}_{q+1}:=R_{\textrm{lin}}+R_{\textrm{cor}}+R_{\textrm{osc}}+R_{\textrm{com}}+R_{\textrm{com}1}.
	\endaligned
\end{equation*}

\subsubsection{Inductive estimate for $\mathring{R}_{q+1}$}
\label{sss:R}
To conclude the proof of Proposition~\ref{p:iteration}, we shall verify the  estimates in \eqref{eq:R} and \eqref{bd:R}. To this end, we estimate each term in the definition of $\mathring{R}_{q+1}$ separately.

In the following we choose $p=\frac{32}{32-7\alpha}>1$ so that it holds  in particular that $r_\perp^{2/p-2}r_\|^{1/p-1}\leq \lambda_{q+1}^\alpha$.
For the linear  error we  obtain
\begin{equation*}
	\aligned
	\|R_{\textrm{lin}}\|_{C_tL^p}
	&\lesssim\|\mathcal{R}\Delta w_{q+1}\|_{C_tL^p}+\|\mathcal{R}\partial_t(w_{q+1}^{(p)}+w_{q+1}^{(c)})\|_{C_tL^p}\\
	&\qquad+\|(v_\ell+z_\ell)\mathring{\otimes}w_{q+1}+w_{q+1}\mathring{\otimes}(v_\ell+z_\ell)\|_{C_tL^p}
	\\
	&\lesssim\|w_{q+1}\|_{C_tW^{1,p}}+\sum_{\xi\in\Lambda}\|\partial_t\textrm{curl}
	(a_{(\xi)}V_{(\xi)})\|_{C_tL^p}\\
	&\qquad+(\|v_q\|_{C_{[t-1,t+1]}L^\infty}+\|z_q\|_{C_{[t-1,t+1]}L^\infty})\|w_{q+1}\|_{C_tL^p},
	\endaligned
\end{equation*}
where by \eqref{bounds} and \eqref{estimate aN}
\begin{equation*}
	\aligned
	\sum_{\xi\in\Lambda}\|\partial_t\textrm{curl}
	(a_{(\xi)}V_{(\xi)})\|_{C_tL^p}&\leq \sum_{\xi\in\Lambda}\left(\|
	a_{(\xi)}\|_{C_tC^1_x}\|\partial_t V_{(\xi)}\|_{C_tW^{1,p}}+\|\partial_ta_{(\xi)}\|_{C_tC^1_x}\| V_{(\xi)}\|_{C_tW^{1,p}}\right)\\
	&\lesssim (\|\mathring{R}_q\|_{C_{[t-1,t+1]}L^1}+1)^{{2}} \ell^{-13}r_{\perp}^{2/p}r_{\|}^{1/p-3/2}\mu\\
	&\qquad+(\|\mathring{R}_q\|_{C_{[t-1,t+1]}L^1}+1)^3 \ell^{-19}r_{\perp}^{2/p-1}r_{\|}^{1/p-1/2}\lambda_{q+1}^{-1}.
	\endaligned
\end{equation*}
In view of \eqref{eq:w1p} as well as  (\ref{corr temporal}),  we deduce
\begin{equation*}
	\aligned
	\|R_{\textrm{lin}}\|_{C_tL^p}
	&\lesssim \left(\lambda_{q+1}^{5\alpha-1/7}+\lambda_{q+1}^{9\alpha-2/7}+\lambda_{q+1}^{27\alpha-1/7}+\lambda_{q+1}^{39\alpha-15/7}\right)(\|\mathring{R}_q\|_{C_{[t-1,t+1]}L^1}+1)^{4}
	\\&+(\|v_q\|_{C_{[t-1,t+1]}L^\infty}+\|z_q\|_{C_{[t-1,t+1]}L^\infty})\lambda_{q+1}^{5\alpha-8/7}(\|\mathring{R}_q\|_{C_{[t-1,t+1]}L^1}+1)^{3}.
	\endaligned
\end{equation*}

The corrector error  is estimated  using \eqref{principle est1}, \eqref{correction est}, \eqref{temporal est1} as
\begin{equation*}
	\aligned
	\|R_{\textrm{cor}}\|_{C_tL^p}
	&\leq\|w_{q+1}^{(c)}+ w_{q+1}^{(t)}\|_{C_tL^{2p}}\| w_{q+1}\|_{C_tL^{2p}}+\|w_{q+1}^{(c)}+ w_{q+1}^{(t)}\|_{C_tL^{2p}}\| w_{q+1}^{(p)}\|_{C_tL^{2p}}
	\\
	&\lesssim \left(\ell^{-19}r_\perp^{1/p}r_\|^{1/(2p)-3/2}
	+\ell^{-4}r_\perp^{1/p-1}r_\|^{1/(2p)-2}\lambda_{q+1}^{-1}\right)\ell^{-2}r_\perp^{1/p-1}r_\|^{1/(2p)-1/2}(\|\mathring{R}_q\|_{C_{[t-1,t+1]}L^1}+1)^{7/2}
	\\&+\|w_{q+1}^{(c)}+ w_{q+1}^{(t)}\|_{C_tL^{2p}}^2
	\\
	&\lesssim \lambda_{q+1}^{-1/7+13\alpha}(\|\mathring{R}_q\|_{C_{[t-1,t+1]}L^1}+1)^{7/2}
	+ \lambda_{q+1}^{-2/7+17\alpha}(\|\mathring{R}_q\|_{C_{[t-1,t+1]}L^1}+1)^{6}.
	\endaligned
\end{equation*}

We continue with the oscillation error $R_{\textrm{osc}}$.
Using \eqref{eR}, \eqref{eB} and the definition of $W_{(\xi)}$ we have
\begin{equation*}
	\aligned
	\|R_{\textrm{osc}}^{(x)}\|_{C_{t}L^p}&\leq \sum_{\xi\in\Lambda}\big\|\mathcal{B}\big(\nabla a^2_{(\xi)},\mathbb{P}_{\geq r_\perp\lambda_{q+1}/2}(W_{(\xi)}\otimes W_{(\xi)})\big)\big\|_{C_{t}L^p}
	\\
	&\lesssim \sum_{\xi\in\Lambda}\|\nabla a^2_{(\xi)}\|_{C_tC^1}\|\cR(W_{(\xi)}\otimes W_{(\xi)})\|_{C_{t}L^p}\lesssim \sum_{\xi\in\Lambda}\|\nabla a^2_{(\xi)}\|_{C_tC^1}\frac{\|W_{(\xi)}\otimes W_{(\xi)}\|_{C_{t}L^p}}{r_\perp\lambda_{q+1}}
	\\&\lesssim \sum_{\xi\in\Lambda}\|\nabla a^2_{(\xi)}\|_{C_tC^1}\frac{\|W_{(\xi)}\|^2_{C_{t}L^{2p}}}{r_\perp\lambda_{q+1}}
	\lesssim \ell^{-{21}}(\|\mathring{R}_q\|_{C_{[t-1,t+1]}L^1}+1)^4r_\perp^{2/p-2}r_\|^{1/p-1}(r_\perp^{-1}\lambda_{q+1}^{-1})\\
	&\lesssim \lambda_{q+1}^{{43}\alpha-1/7}(\|\mathring{R}_q\|_{C_{[t-1,t+1]}L^1}+1)^4.
	\endaligned
\end{equation*}
For the second term $R_{\textrm{osc}}^{(t)}$ we use Fubini's theorem to integrate along the orthogonal directions of $\phi_{(\xi)}$ and $\psi_{(\xi)}$ and apply (\ref{bounds}) to deduce
\begin{equation*}
	\aligned
	\|R_{\textrm{osc}}^{(t)}\|_{C_{t}L^p}&\leq \mu^{-1}\sum_{\xi\in\Lambda}\|\partial_t a_{(\xi)}^2\|_{C^{0}_{t,x}}\|\phi_{(\xi)}\|_{C_{t}L^{2p}}^2\|\psi_{(\xi)}\|_{C_{t}L^{2p}}^2
	\\
	&\lesssim (\|\mathring{R}_q\|_{C_{[t-1,t+1]}L^1}+1)^{5/2} \mu^{-1}\ell^{-15}r_\perp^{2/p-2}r_\|^{1/p-1}\lesssim \lambda_{q+1}^{31\alpha-9/7}(\|\mathring{R}_q\|_{C_{[t-1,t+1]}L^1}+1)^{5/2}.
	\endaligned
\end{equation*}

In view of the standard mollification estimates and \eqref{inductionv} it holds
\begin{equation*}
	\aligned
	\|R_{\textrm{com}}\|_{C_{t}L^1}&\lesssim \ell\|v_q\|_{C^1_{[t-1,t+1],x}}(\|v_q\|_{C_{[t-1,t+1]}L^2}+\|z_q\|_{C_{[t-1,t+1]}L^2})\\
	&\quad+\ell^{{1}/{2}-\delta}(\|z_q\|_{C_{[t-1,t+1]}^{{1}/{2}-\delta}L^2}+\|z_q\|_{C_{[t-1,t+1]}H^{1-\delta}})(\|v_q\|_{C_{[t-1,t+1]}L^2}
	+\|z_q\|_{C_{[t-1,t+1]}L^2}),
	\endaligned
\end{equation*}
where $\delta<\frac{1}{12}$.
Finally, we use Proposition \ref{fe z} and \eqref{eR} to obtain
\begin{align*}\|\mathcal{R}(z_\ell-z_{q+1})\|_{C_tL^2}\lesssim\|z_\ell-z_{q+1}\|_{C_tL^2}&\leq \|z_\ell-z_q\|_{C_tL^2}+\|z_q-z_{q+1}\|_{C_tL^2}
\\&\leq\ell^{\frac{1}{2}-\delta}\|z\|_{C_{[t-1,t+1]}^{1/2-\delta}L^2}+\|z\|_{C_{t}H^{1-\delta}}\lambda_{q+1}^{-\frac\alpha8(1-\delta)}
\end{align*}
which implies
{\begin{equation*}
		\aligned
		\|R_{\textrm{com}1}\|_{C_tL^1}&\lesssim (\|v_{q+1}\|_{C_tL^2}+\|z_{q+1}\|_{C_{[t-1,t+1]}L^2}
		+\|z_q\|_{C_{[t-1,t+1]}L^2}+1)\|z_\ell-z_{q+1}\|_{C_tL^2}
		\\&\leq  (\ell^{\frac{1}{2}-\delta}\|z\|_{C_{[t-1,t+1]}^{1/2-\delta}L^2}+\|z\|_{C_{[t-1,t+1]}H^{1-\delta}}\lambda_{q+1}^{-\frac\alpha8(1-\delta)})(\|v_{q+1}\|_{C_tL^2}+\|z\|_{C_{[t-1,t+1]}L^2}
		+1)
		\\&\leq  \$z\$\lambda_{q+1}^{-\frac\alpha8(1-\delta)}(\|v_{q+1}\|_{C_tL^2}+\|z\|_{C_{[t-1,t+1]}L^2}
		+1).
		\endaligned
\end{equation*}}
Here $\$z\$=\|z\|_{C_{[t-1,t+1]}^{1/2-\delta}L^2}+\|z\|_{C_{[t-1,t+1]}H^{1-\delta}}$.
Summing up all the above estimates, we obtain
\begin{equation*}
	\aligned
	\|\mathring{R}_{q+1}\|_{C_{t}L^1}&\lesssim \lambda_{q+1}^{{43}\alpha-1/7}(\|\mathring{R}_q\|_{C_{[t-1,t+1]}L^1}+1)^{4}+ \lambda_{q+1}^{-2/7+17\alpha}(\|\mathring{R}_q\|_{C_{[t-1,t+1]}L^1}+1)^{6}
	\\&\quad+(\|v_q\|_{C_{[t-1,t+1]}L^\infty}+\|z_q\|_{C_{[t-1,t+1]}L^\infty})\lambda_{q+1}^{5\alpha-8/7}(\|\mathring{R}_q\|_{C_{[t-1,t+1]}L^1}+1)^{3}
	\\&\quad+\$z\$\lambda_{q+1}^{-\frac\alpha8(1-\delta)}(\|v_{q+1}\|_{C_tL^2}+\|z\|_{C_{[t-1,t+1]}L^2}+1)
	\\&\quad+(\ell\|v_q\|_{C^1_{[t-1,t+1],x}}+\ell^{\frac{1}{2}-\delta}\$z\$)(\|v_q\|_{C_{[t-1,t+1]}L^2}+\|z\|_{C_{[t-1,t+1]}L^2}).
	\endaligned
\end{equation*}
Thus taking the $r$-th moment, using H\"older's inequality and  \eqref{vqc1}, \eqref{bd:R}, \eqref{z ps} we obtain
\begin{equation*}
	\aligned
	\$\mathring{R}_{q+1}\$_{L^1,r}&\lesssim \lambda_{q+1}^{{43}\alpha-1/7}(\$\mathring{R}_q\$^4_{L^1,4r}+1)+ \lambda_{q+1}^{-2/7+17\alpha}(\$\mathring{R}_q\$^{6}_{L^1,6r}+1)
	\\&\quad+(\$v_q\$_{C_{t,x}^1,2r}+\$z_q\$_{L^\infty,2r})\lambda_{q+1}^{5\alpha-8/7}(\$\mathring{R}_q\$^3_{L^1,6r}+1)
	\\&\quad+\lambda_{q+1}^{-\frac\alpha8(1-\delta)}(\$v_{q+1}\$_{L^2,2r}+(2r-1)^{1/2}L)L(2r-1)^{1/2}
	\\&\quad+(\ell\$v_q\$_{C^1_{t,x},2r}+\ell^{\frac{1}{2}-\delta}L(2r-1)^{1/2})(\$v_q\$_{L^2,2r}+L(2r-1)^{1/2})
	\\&\lesssim\lambda_{q+1}^{{43}\alpha-1/7} {\lambda^{12}_q}(6^q\cdot 24L^2r)^{6^{q+1}}+\lambda_{q+1}^{5\alpha-8/7}(6^q\cdot 24L^2r)^{3 (6^q)}{\lambda^6_q}(\lambda_q^4+\lambda_{q+1}^{\alpha/8}L)
	\\&\quad+(M_0\bar{e}^{1/2}+L)(L\lambda_{q+1}^{-\frac\alpha8(1-\delta)}+\ell\lambda_q^{4}+\ell^{\frac{1}{2}-\delta} L)
	\\&\lesssim\lambda_{q+1}^{{43}\alpha-1/7} \lambda_q^{{18}}+\lambda_{q+1}^{5\alpha-8/7}\lambda_q^{{9}}(\lambda_q^{4}+\lambda_{q+1}^{\alpha/8})+\lambda_{q+1}^{-\frac\alpha8(1-\delta)}
	\\&\leq\frac{1}{48}\delta_{q+3}\underline{e}.
	\endaligned
\end{equation*}
Here in the third inequality we used $(6^q\cdot 24L^2r)^{6^{q}}\leq \lambda_q$ and
in the last inequality we used
${43}\alpha<1/14$ and ${18}/b+2\beta b^2<1/14$ and $\alpha>40\beta b^2$ and $\alpha b>32/7$. Hence \eqref{eq:R} holds on the level $q+1$.

Similarly, for $m\geq1$ we use the first inequality as above with $r$ replaced by $m$ and instead of \eqref{inductionv} we use \eqref{inductionv C1} and \eqref{inductionv m} to obtain
\begin{equation*}
	\aligned
	\$\mathring{R}_{q+1}\$_{L^1,m}&\lesssim \lambda_{q+1}^{{43}\alpha-1/7}(\$\mathring{R}_q\$^4_{L^1,4m}+1)+ \lambda_{q+1}^{-2/7+17\alpha}(\$\mathring{R}_q\$^{6}_{L^1,6m}+1)
	\\&\quad+(\$v_q\$_{C_{t,x}^1,2m}+\$z_q\$_{L^\infty,2m})\lambda_{q+1}^{5\alpha-8/7}(\$\mathring{R}_q\$^3_{L^1,6m}+1)
	\\&\quad+\lambda_{q+1}^{-\frac\alpha8(1-\delta)}(\$v_{q+1}\$_{L^2,2m}+(2m)^{1/2}L)(2m)^{1/2}L
	\\&\quad+\Big(\ell\$v_q\$_{C^1_{t,x},2m}+\ell^{\frac{1}{2}-\delta}L(2m)^{1/2}\Big)(\$v_q\$_{L^2,2m}+L(2m)^{1/2})
	\\&\lesssim\lambda_{q+1}^{{43}\alpha-1/7} (6^q\cdot 24mL^2)^{6^{q+1}}{\lambda^{12}_q}\\
	&\quad+\lambda_{q+1}^{5\alpha-8/7}(6^q\cdot 24mL^2)^{3 (6^q)}{\lambda^6_q}\Big(\lambda_q^{23/7}(6^{q-1}\cdot 32mL^2)^{4 (6^{q-1})}+\lambda_{q+1}^{\alpha/8}(2m)^{1/2}L\Big)
	\\&\quad+(2m)^{1/2}L\lambda_{q+1}^{-\frac\alpha8(1-\delta)}\Big(({6^q\cdot 24mL^2)^{3 (6^q)}}{\lambda_{q+1}}+\bar{e}^{1/2}+(2m)^{1/2}L\Big)
	\\&\quad+\Big({(6^q\cdot 4mL^2)^{3 (6^{q-1})}}{\lambda_q}+\bar{e}^{1/2}+(2m)^{1/2}L\Big)
	\\&\quad\times\Big(\ell\lambda_q^{23/7}(6^{q-1}\cdot 32mL^2)^{4 (6^{q-1})}+\ell^{\frac{1}{2}-\delta} L(2m)^{1/2}\Big)
	\\&\leq(6^{q+1}\cdot 4mL^2)^{6^{q+1}}{\lambda_{q+1}^2}.
	\endaligned
\end{equation*}

The proof of Proposition~\ref{p:iteration} is therefore complete.

\section{Stationary solutions to the stochastic Navier--Stokes system}
\label{s:4}

We recall that  the trajectory space is  $\mathcal{T}= C(\mR;L^2_{\sigma})\times C(\mR;L_{\sigma}^2)$ and the corresponding shifts $S_t$, $t\in\mR$,  on trajectories are given by
$$
S_t(u,B)(\cdot)=(u(\cdot+t),B(\cdot+t)-B(t)),\quad t\in\mR,\quad (u,B)\in\cT.
$$
The notion of stationary solution was introduced in Definition~\ref{d:1.1}. Our first result of this section is existence of stationary solutions as limits of ergodic averages of solutions constructed in the previous section. This in particular implies their non-uniqueness.

\bt\label{th:s1}
Let $u$ be a solution obtained in Theorem~\ref{thm:6.1} with $e(t)=K$ for some $K\geq 8\cdot 48 L^2r$ and satisfying \eqref{est:u1} and \eqref{es:v1} with given $\eps>0$ and $r>1$. Then there exists a sequence $T_{n}\to\infty$ and
a stationary  solution $((\tilde\Omega,\tilde{\mathcal{F}},\tilde{\mathbf{P}}),\tilde u,\tilde B)$ to \eqref{1} such that
$$
\frac{1}{T_{n}}\int_{0}^{T_{n}}\mathcal{L}[S_{t}(u,B)] \dif t\to \mathcal{L}[\tilde u,\tilde B]
$$
weakly in the sense of probability measures on $\mathcal{T}$ as $n\to\infty$. Moreover, it holds true that
\begin{align}\label{eq:s}
	\begin{aligned}
		\tilde{\mathbf{E}}\|\tilde u\|_{L^2}^2=K,
	\end{aligned}
\end{align}
and for $\eps>0$
\begin{align}\label{eq:s1}
	\begin{aligned}
		\$\tilde u-\tilde z\$_{W^{1,1},r}\leq \eps,
	\end{aligned}
\end{align}
for $\tilde z(t)=\mathbb{P}\int_{-\infty}^te^{(t-s)(\Delta-1)}\dif \tilde B_s$ and for some $\vartheta>0$ and for every $N\in\mathbb{N}$
\begin{align}\label{bd:u}
	\tilde{\E}\sup_{t\in[-N,N]}\|\tilde u(t)\|^{2{r}}_{H^\vartheta}+\E\|\tilde u\|^{2{r}}_{C^\vartheta([-N,N], L^2)}\lesssim N.		
\end{align}
\et

\begin{proof}
	We define the ergodic averages of the solution $(u,B)$ as the probability measures on the space of trajectories $\cT$
	\begin{align*}
		\nu_T=\frac1T\int_0^T\mathcal{L}[S_t(u,B)]\dif t,\qquad T\geq 0.
	\end{align*}
	By Theorem \ref{thm:6.1} we obtain
	\begin{align*}
		\sup_{s\in\mR}\E\sup_{t\in[-N,N]}\|u(t+s)\|^{2r}_{H^\vartheta}&\leq\sup_{s\in\mR}\sum_{i=-N}^{N-1}\E\sup_{t\in[i,i+1]}\|u(t+s)\|_{H^\vartheta}^{2r}\no
		\\&\leq 2N\sup_{s\geq0}\E\sup_{t\in[0,1]}\|u(t+s)\|^{2r}_{H^\vartheta}
		\lesssim N,
	\end{align*}
	and similarly
	\begin{align*}
		\sup_{s\in\mR}\E\|u(\cdot+s)\|^{2}_{C^\vartheta([-N,N], L^2)}\lesssim N.
	\end{align*}
	For $R_{N}>0$, $N\in\mathbb{N}$, we note that the set
	\begin{align*}
		K_M:=\cap_{N=M}^\infty\bigg\{g_1;\,&\|g_1\|_{C^\vartheta_{[-N,N]} L^2}+\sup_{t\in[-N,N]}\|g_1(t)\|_{H^\vartheta}\leq R_N\bigg\}
	\end{align*}
	is relatively compact in $C(\mR;{L^2_\sigma})$. As a consequence, we deduce that the time shifts $S_tu$, $t\in\mR$,  are tight on $ C(\mR;{L^2_\sigma})$. Since  $S_tB$ is a Wiener process for every $t\in\mR$, the law of  $S_tB$ does not change with $t\in\mR$ and is tight. Accordingly, for any $\eps>0$ there is a compact set $\bar K_\eps$ in $\cT$ such that
	$$\sup_{t\in\mR}\bP(S_{t} (u,B) \in \bar{K}^c_\eps)<\eps.$$
	This implies
	\begin{align*}
		\nu_T (\bar{K}_\eps^c) =& \frac{1}{T} \int_0^T \mathbf{P} (S_{t} (u,B) \in \bar{K}^c_\eps) \dif t
		< \eps
	\end{align*}
	and therefore there is a weakly converging subsequence of the probability measures $\nu_{T}$, $T\geq0$. That is, there is a subsequence $T_n\to\infty$ and $\nu\in \cP(\cT)$ such that $\nu_{T_n}\to \nu$ weakly in $\cP(\cT)$.

	Define  the set
	\begin{align*}
		A&=\Big\{(u,B)\in \cT;\
		\langle u(t),\psi\rangle+\int_s^t\<\div(u\otimes u), \psi\>\dif r
		\\&\qquad\qquad= \langle u(s),\psi\rangle + \int_{s}^{t}\langle \Delta u
		,\psi\rangle\dif r+\<B(t)-B(s),\psi\>,
		\quad\forall \psi \in C^{\infty}(\mathbb{T}^{3}),\ \div\psi=0,\  t\geq s\Big\}.
	\end{align*}
	Since $(u,B)$ in the statement of the theorem satisfies the equation, we
	have for all $t\in\mR$
	\begin{align*}
		\cL[S_{t}(u,B)](A)=1.
	\end{align*}
	Hence, also $\nu_{T_{n}}(A)=1$ for all $n\in\mathbb{N}$. By Jakubowski--Skorokhod representation theorem, there is a probability space $(\tilde \Omega,\tilde{\mathcal{F}},\tilde{\mathbf{P}})$ with  a sequence of random variables $(\tilde u^n,\tilde B^n)$, $n\in\mathbb{N}$, such that $\cL[\tilde u^n,\tilde B^n] ={\nu_{T_{n}}}$ %
	and  $(\tilde u^n,\tilde B^n)$ satisfy equation \eqref{1} on $\mR$. By \eqref{eq:K1} we know
	\begin{align}\label{eq:un1}
		\tilde{\E}\|\tilde{u}^n(t)\|_{L^2}^2=\frac1{T_n}\int_0^{T_n}\E\|{S_su(t)}\|_{L^2}^2\dif s=K.
	\end{align}
	By \eqref{est:u1}
	\begin{align}\label{eq:un2}
		\sup_n\tilde \E\|\tilde{u}^n(t)\|_{L^2}^{2r}=\sup_n\frac1{T_n}\int_0^{T_n}\E\|S_su(t)\|_{L^2}^{2r}\dif s<\infty.
	\end{align}
	Moreover, there is a  random variable $(\tilde u,\tilde B)$ having the law $\cL[\tilde u,\tilde B]=\nu$ so that
	$$
	(\tilde u^n,\tilde B^n)\to(\tilde u,\tilde B)\qquad \tilde{\mathbf{P}}\text{-a.s. in }\mathcal{T}.
	$$
	Thus, we can pass to the limit in the equation to deduce that $\nu$ is a law of a solution on $\mR$.
	
	The shift invariance  follows from the same argument as in \cite[Lemma 5.2]{BFH20e}. Namely, it holds for $G\in C_b(\cT)$ and $r\in\mR$
	\begin{align*}
		\int_{\cT}G\circ S_r(u,B)\dif \nu(u,B)=\int_{\cT}G (u,B)\dif \nu(u,B).
	\end{align*}
	Finally,  \eqref{eq:s} follows from \eqref{eq:un1}, \eqref{eq:un2} and \eqref{eq:s1}, \eqref{bd:u} follow from a lower-semicontinuity argument and the related bound for the approximations. In fact, we define
$$\tilde{z}^n(t):=\int_{-\infty}^te^{(t-s)(\Delta-I)}\dif \tilde B^n=\tilde B^n(t)+\int_{-\infty}^t(\Delta-I)e^{(t-s)(\Delta-I)}\tilde B^n\dif s.$$
We know that  $\tilde{z}(t)^n\rightarrow\tilde{z}(t)=\int_{-\infty}^te^{(t-s)(\Delta-I)}\dif \tilde B$ in $C(\mathbb{R},H^{-2})$ $\tilde{\mathbf{P}}$ a.s. Thus \eqref{eq:s1} follows from  lower-semicontinuity.
\end{proof}

Using the above result and choosing different $K$, the first claim in Theorem \ref{th:main} follows.

By a general result applied also in \cite{BFH20e,FFH21,HZZ22} and using Theorem \ref{th:s1} we obtain existence of infinitely many ergodic stationary solutions as follows.

\bt\label{thm:4.5}
Let $r>1$. For $K\geq 8\cdot 48 L^2r$ there exists $C>0$ and  an ergodic stationary solution $((\Omega,\mathcal{F},\mathbf{P}),u,B)$ satisfying
\begin{align}\label{eq:s55}
	\begin{aligned}
		\mathbf{E}\| u\|_{L^2}^2= K,
	\end{aligned}
\end{align}
and for some $\vartheta>0$ and for every $N\in\mathbb{N}$
\begin{align}\label{eq:s56}
\E\sup_{t\in[-N,N]}\| u(t)\|^{{2r}}_{H^\vartheta}+\E\| u\|^{{2r}}_{C^\vartheta([-N,N], L^2)}\leq C N.
\end{align}
\et

\begin{proof}
In view of Theorem~\ref{th:s1}, this is a consequence of a Krein--Milman argument. In particular, we observe that the set of all laws of  stationary solutions satisfying \eqref{eq:s55} and \eqref{eq:s56} is non-empty, convex, tight and closed which follows from the arguments in the proof of Theorem~\ref{th:s1}. Hence there exist an extremal point. By a classical contradiction argument,  it is the law of an ergodic stationary solution.

Non-uniqueness of ergodic stationary solutions follows from choosing different $K$.
\end{proof}

\section{Stationary solutions to the stochastic Euler equations}
\label{s:5}

In this section we proceed with a construction of stationary solutions to the stochastic Euler equations \eqref{2}. Furthermore, we prove that the law of every stationary solution to the stochastic Euler equations \eqref{2} is a limit of the law of stationary solutions to the stochastic Navier--Stokes equations with vanishing viscosities. 

In order to  construct stationary solutions to the stochastic Euler equations \eqref{2},
we decompose its solution $u$ as $v+z$ where $z$ solves the linear stochastic problem
\begin{equation}\label{linear E}
	\aligned
	\dif z + z \,\dif t&=\dif B
	\endaligned
\end{equation}
and $v$ is a solution to the nonlinear  equation with random coefficients
\begin{equation}\label{nonlinear E}
	\aligned
	\partial_tv -z+\div((v+z)\otimes (v+z))+\nabla P&=0,
	\\
	\div v&=0.
	\endaligned
\end{equation}
Here, unlike in Proposition~\ref{fe z}, we cannot  use the smoothing effect of the Laplacian and hence we obtain spatial regularity of $z$ from the strengthened assumption on the Wiener process. Specifically, we assume  $\tr((-\Delta)^\sigma GG^*)<\infty$ for some $\sigma>0$. Let  $z$ be the unique stationary solution to \eqref{linear E} for which  the bound from Proposition~\ref{fe z} replaced by: for any $\delta\in(0,1/2)$, $p\geq1$
\begin{align}\label{b:z}
	\sup_{t\in\mR}\E\Big[\|z\|_{C_t^{1/2-\delta}H^\sigma}^p\Big]\leq (p-1)^{p/2}L^p.
\end{align}

\bt\label{th:5.1}
Assume that $\tr((-\Delta)^\sigma GG^*)<\infty$ for some $\sigma>0$. There exist infinitely many stationary solutions $((\Omega,\cF,\bP), u, B)$ to stochastic Euler equations \eqref{2} on $\mR\times \mT^3$. In particular, let $r>1$ and for a given $K\geq8\cdot 48L^2r$  with $L$ being the bound for the noise   in Proposition \ref{fe z}, there exists a stationary solution  $((\Omega,\cF,\bP), u, B)$ to \eqref{2} satisfying
$$
\E \|u\|_{L^2}^2=K,%
$$
as well as \eqref{bd:u} for some $\vartheta>0$.

Moreover, for an arbitrary sequence of vanishing viscosities $\nu_n\to0$, $n\in\N$, there exist a sequence of stationary solutions $u_n$, $n\in\N$, to the following stochastic Navier--Stokes equations
\begin{align}\label{eq:unn}
\dif u_n+\div(u_n\otimes u_n)\,\dif t+\nabla P_{n}\,\dif t=\nu_n\Delta u_n\,\dif t+\dif B,
\end{align}
so that the corresponding family of laws $\mathcal{L}[u_{n}]$, $n\in\N$, is tight in $C(\mR;L_{\sigma}^2)$ and
every accumulation point is a stationary solution to  \eqref{2}.

Finally, there exist infinitely many  ergodic stationary solutions to \eqref{2}.
\et

\begin{proof}
We revisit the convex integration scheme from Section \ref{s:in}. Since the statement of Proposition~\ref{fe z} is now replaced by \eqref{b:z}, the bounds in \eqref{z ps} are changed to
\begin{align}\label{eq:znew}
	\$z_q\$_{L^\infty,p}\leq\$z_q\$_{C_t^{1/2-2\delta}L^\infty,p}\leq \lambda_{q+1}^{\alpha/4}(p-1)^{1/2}L.
\end{align}
In Section \ref{s:en}, we need to modify the bounds \eqref{eq:00} and \eqref{eq:01}. The former one now reads for every $\eps>0$ as
\begin{align*}
	2\mathbf{E}\|(v_\ell +z_{q+1}) \cdot w_{q+1}^{(p)}\|_{L^1}
	&\lesssim (\$v_\ell\$_{L^\infty,2}+\$z_q\$_{L^\infty,2})\$  w_{q+1}^{(p)}\$_{L^1,2}+{\$z_{q+1}-z_q\$_{L^2,2}\$  w_{q+1}^{(p)}\$_{L^{2},2}}
	\\&\lesssim (\lambda_q^{4}+\lambda_{q+1}^{\alpha/4}L)\ell^{-2}\delta_{q+2}^{1/2}r_{\perp}^{1-\eps}r_{\|}^{\frac12(1-\eps)}+{\lambda_{q+1}^{-{\alpha\sigma}/8}M_0\bar{e}^{1/2}}\\
	&\lesssim L\lambda_{q+1}^{5\alpha-\frac87(1-\eps)}+{\lambda_{q+1}^{-{\alpha\sigma}/8}M_0\bar{e}^{1/2}}\leq \frac{1}{96}
	\lambda_{q+1}^{-2\beta b}\leq \frac{\delta_{q+2}}{96}e(t),
\end{align*}
where we use \eqref{estimate wqp} and \eqref{estimate wq1} to control $\$  w_{q+1}^{(p)}\$_{L^{2},2}$ and we need $\alpha \sigma>16\beta b$.
The latter one now relies on
\begin{align*}
	\E\|z_{q+1}-z_q\|_{L^2}^2\leq L\lambda_{q+1}^{-\alpha\sigma/8},
\end{align*}
which requires
$M_0{(\bar{e}^{1/2}+L)}\leq \lambda_{q+1}^{{\alpha\sigma}/{8}-2\beta b}$, i.e. $\alpha\sigma>16\beta b$ and $a$ large enough to absorb the extra constant.

For the control of $\mathring{R}_{q+1}$  we can use \eqref{eq:znew} to derive the same bounds for most of the terms as in Section~\ref{sss:R}. The main change comes from the following two parts in $R_{\text{com}}$ and $R_{\text{com}1}$, namely,
\begin{align*}
	\ell^\sigma \|z_q\|_{C_{[t-1,t+1]}H^\sigma}(\|v_q\|_{C_{[t-1,t+1]}L^2}+\|z_q\|_{C_{[t-1,t+1]}L^2}+{\|v_{q+1}\|_{C_{[t-1,t+1]}L^2}}),
\end{align*}
and
\begin{align*}
	\|z\|_{C_{[t-1,t+1]}H^\sigma}\lambda_{q+1}^{-\alpha\sigma/8}(\|v_q\|_{C_{[t-1,t+1]}L^2}+\|z_q\|_{C_{[t-1,t+1]}L^2}+1).
\end{align*}
Then, when we estimate $\$\mathring{R}_{q+1}\$_{L^1,r}$ we
have the following extra term
\begin{align*}
	{(M_0\bar e^{1/2}+L)L(\lambda_{q + 1}^{-\alpha\sigma/8}+\ell^\sigma),}
\end{align*}
which requires ${(M_0\bar e^{1/2}+L)}L\leq \lambda_{q+1}^{{\alpha\sigma}/8-2b\beta^2}$, i.e. $\alpha \sigma >16\beta b^2$. We obtain this additional  estimate by choosing $\beta$ small enough.

Consequently, we deduce the existence and non-uniqueness of  solutions to the stochastic Euler equations as in Theorem \ref{thm:6.1}. Furthermore, the existence and non-uniqueness of stationary solutions  follow by the same argument as in Theorem \ref{th:s1}. This completes the proof of the first claim in Theorem~\ref{th:5.1}.

For the second result in Theorem~\ref{th:5.1}, we first apply Theorem \ref{th:s1} to derive the existence of stationary solutions $u_n$, $n\in\N$, to equations \eqref{eq:unn} with $u_n$ satisfying \eqref{bd:u} uniformly in $n$. More precisely, we choose  $z_n$ satisfying
\begin{equation}\label{linear En}
	\aligned
	\dif z_n + z_n \,\dif t&=\nu_n \Delta z_n\,\dif t+ \dif B.
	\endaligned
\end{equation}
Then $z_n$ satisfies \eqref{b:z} uniformly in $n$, using again the regularity of the noise instead of the smoothing effect of the Laplacian. Hence, by exactly the same argument as above and Theorem \ref{th:s1} we know that \eqref{bd:u} holds uniformly in $n$, which implies  tightness of $u_n$, $n\in\N$, in $C(\mR;L_{\sigma}^2)$. By Jakubowski--Skorokhod representation theorem, we can modify the stochastic basis and pass to the limit in the approximate Navier--Stokes equations \eqref{eq:unn} and derive the claim.

Existence and non-uniqueness of ergodic stationary solutions follows from the same argument as in Theorem~\ref{thm:4.5}.
\end{proof}

In the following, we show that for every solution to the stochastic Euler equation satisfying suitable moment bounds we can find a solution to stochastic Navier--Stokes equations close to it.  The proof follows from essentially the same argument as \cite[Theorem 1.3]{BV19a}. However, in the stochastic setting, it is necessary to control any $m$th moment of the solutions during the iteration in Proposition \ref{p:iteration}, even though the given  solution to the Euler equations is not  required to possess finite moments of all orders.  Accordingly, we cut-off the solution to the Euler equations in a suitable way to have all finite moments and start the iteration  from this truncated solution.

\bt\label{thm:new}Assume that $\tr((-\Delta)^\sigma GG^*)<\infty$ for some $\sigma>0$.
Let $r>1$, $\eps>0$ be fixed and $u$ be an analytically weak solution to the stochastic Euler equations \eqref{2}  satisfying for some $\bar\beta'>0$, $r'>r$ and $L \geq (2\pi)^3$
\begin{equation}\label{eq:new}
\$u\$_{H^{\bar\beta'},2r'}+\$u\$_{C_t^{\bar \beta'}L^2,2r'}\leq L.
\end{equation}
There exist $\nu>0$, $\vartheta>0$ and an analytically weak solution $u_\nu$ to the stochastic Navier--Stokes equations \eqref{1}  on $\mR\times \mT^3$ which belongs to $ {C(\mR;H^{\vartheta})\cap C^{\vartheta}(\mR;L^{2})}$ $\mathbf{P}$-a.s.,
$$
\$u_\nu\$_{H^\vartheta,2r}+\$u_\nu\$_{C_{t}^\vartheta L^2,2r}<\infty,
$$
and
\begin{align}\label{eq:uwp new}
\$u_\nu - u\$_{L^2,r}\leq \varepsilon.
\end{align}
\et

\begin{proof}  Choose all the parameters except $a,\beta$ as in Proposition \ref{p:iteration}. Without loss of generality we assume that $\bar\beta'\leq 3/4$. We take $2\beta<\min(\frac{\bar\beta\wedge\sigma\wedge\bar\beta'}{b^2},\frac{\sigma\alpha}{8b}) $ with
	$$\bar\beta<\frac{4\bar\beta'}{3}\wedge\Big(\frac{r'}{r}-1\Big)\wedge \Big(\frac{2(\bar\beta'\wedge\sigma)}{3-2(\bar \beta'\wedge\sigma)}\Big),$$ and $n$ large enough such that in \eqref{eqtheta} for $a\geq a_{0}>0$
\begin{align}\label{eqtheta1}C\sum_{q\geq n}\delta_{q+1}^{\frac{1-\vartheta}{2}}\lambda_{q+1}^{4\vartheta}\leq \eps/2,
\end{align}
where $C$ is the implicit constant in \eqref{eqtheta}. We split $u=v+z$ where $z$ solves \eqref{linear E}
whereas $v$ is a solution to the nonlinear equation \eqref{nonlinear E}. We first have by \eqref{eq:new} and \eqref{b:z}
\begin{align}\label{bdvmom}\$v\$_{H^{\bar\beta'\wedge\sigma},2r'}+\$v\$_{C_t^{\bar \beta'\wedge(\frac12-\delta)}L^2,2r'}\leq CL,\end{align}
for $\delta\in (0,1/2)$. Here we could take the larger $L$ in  \eqref{eq:new} and \eqref{b:z}.
Define
$$z_n=\mathbb{P}_{\leq \lambda_{n+1}^{\alpha/8}}z.$$
Recall the bound for $z$ from \eqref{b:z}. Hence, we have
\begin{align}\label{bd:s5z1}
	\$z_n-z\$_{L^2,2r}\lesssim \lambda_{n+1}^{- \alpha\sigma/8}.
\end{align}
 We
then define
$$v^{(\lambda_n)}:=\chi_{\lambda_n}(v),$$
with $\chi_{\lambda_n}\in C^1(\mR^3)$ and $\chi_{\lambda_n}(x)=x$ for $|x|\leq \lambda_n$, $|\chi_{\lambda_n}|\leq {\lambda_n}+1$, {$\chi_{\lambda_n}(x)\leq |x|$ for $x\in \mR^3$,} and $|\chi_{\lambda_n}'|\leq 1$,  $\supp\chi_{\lambda_n}'\subset \{|x|\leq {\lambda_n}+1\}$.
We then obtain
for $q=\frac3{2\bar\beta'}\geq 2$ and $\bar r=\frac{r'r}{r'-r}$
\begin{align}\label{bd:vM2}
	\$\chi_{\lambda_n}'(v)-1\$_{L^q,\bar r}&\lesssim \Big(\E\big(\int_{|v|>\lambda_n}\dif x\big)^{\bar r/q}\Big)^{1/\bar r}\lesssim \lambda_n^{-\bar\beta}\Big(\E\big(\int|v|^{\bar \beta q}\dif x\big)^{\bar r/q}\Big)^{1/\bar r}
\\&\lesssim  \lambda_n^{-\bar\beta}\Big(\E\|v\|_{L^2}^{\bar \beta\bar r}\Big)^{1/\bar r}\lesssim\lambda_n^{-\bar \beta}\$v\$_{H^{\bar\beta'\wedge\sigma},2r'}^{\bar \beta},\no
\end{align}
where in the third inequality we used H\"{o}lder's inequality and $\bar \beta<\frac{4\bar\beta'}{3}$, and we used $\bar \beta< 2(\frac{r'}{r}-1)$ to have $\bar r\bar \beta/q\leq 2r' $ in the last inequality.
Moreover, we have
\begin{align}\label{bd:vM1}
	\$ v^{({\lambda_n})}-v\$_{L^2,2r}&\lesssim \Big(\E\big(\int_{|v|>\lambda_n}|v|^2\dif x\big)^{r}\Big)^{\frac1{2r}}\lesssim \lambda_n^{-\bar\beta}\Big(\E\big(\int|v|^{2\bar \beta+2}\dif x\big)^{r}\Big)^{\frac1{2r}}
\\&\lesssim \lambda_n^{-\bar\beta}\Big(\E\|v\|_{H^{\bar\beta'\wedge\sigma}}^{r(2\bar \beta+2)}\Big)^{\frac1{2r}}\lesssim\lambda_n^{-\bar \beta}\$v\$_{H^{\bar\beta'\wedge\sigma},2r'}^{\frac12(\bar\beta+1)},\no
\end{align}
where in the last inequality we used Sobolev embedding   and $\bar \beta<(\frac{r'}{r}-1)\wedge (\frac{2(\bar\beta'\wedge\sigma)}{3-2(\bar \beta'\wedge\sigma)})$ to have $r(2\bar \beta+2)<2r'$ and $\|v\|_{L^{2\bar \beta+2}}\lesssim\|v\|_{H^{\bar\beta'\wedge\sigma}}$.
Then we
define
$$v_n=(v^{(\lambda_n)}*_x\phi_{\lambda_{n}^{-1}})*_t\psi_{\lambda_{n}^{-1}}$$
where $\phi$ and $\psi$ were introduced in Section~\ref{s:p}. We
 apply  \eqref{bd:s5z1}, \eqref{bd:vM1}  to bound
\begin{align}\label{eq:app}
	\$ u-v_n-z_n\$_{L^2,2r}&\leq \$ z-z_n\$_{L^2,2r}+\$v-v^{(\lambda_{n})}\$_{L^2,2r}+\$v^{(\lambda_{n})}-v_n\$_{L^2,2r}
\\&\lesssim\lambda_{n+1}^{- \alpha\sigma/8} +\lambda_n^{-\bar \beta}+\lambda_n^{-(\bar\beta'\wedge\sigma\wedge (\frac12-\delta))}\leq \frac{\eps}2.\no
\end{align}
Here in the second inequality we used
\begin{align}\label{bdvn}
\$v^{(\lambda_{n})}-v_n\$_{L^2,2r}&\lesssim 	\lambda_n^{-(\bar\beta'\wedge\sigma)}\$v^{(\lambda_{n})}\$_{H^{\bar\beta'\wedge\sigma},2r'}+\lambda_n^{-(\bar\beta'\wedge (\frac12-\delta))}\$v^{(\lambda_{n})}\$_{C_t^{\bar \beta'\wedge(\frac12-\delta)}L^2,2r'}\no
\\&\lesssim\lambda_n^{-(\bar\beta'\wedge\sigma)}\$v\$_{H^{\bar\beta'\wedge\sigma},2r'} +\lambda_n^{-(\bar\beta'\wedge (\frac12-\delta))}\$v\$_{C_t^{\bar \beta'\wedge(\frac12-\delta)}L^2,2r'}\lesssim \lambda_n^{-(\bar\beta'\wedge\sigma\wedge (\frac12-\delta))},
	\end{align}  and in the last step we choose $a$ large enough and also $n$ large enough.
Also define the energy function
$$e(t):=\frac1{1-\delta_{n+1}}\E\|v_n(t)+z_n(t)\|_{L^2}^2.$$
We fix such $n$ and start  the iteration with $(v_n,z_n)$ at step $n$ and  \eqref{induction ps} on the level $q=n$ reads as
\[ \partial_t v_n-z_n + \tmop{div} ((v_n+z_n) \otimes (v_n+z_n)) + \nabla p_n- \lambda_{n}^{-2}\Delta (v_n+z_n)=: \tmop{div} \mathring{R}_n \]
with
\begin{align*}
\mathring{R}_n &= (v_n+z_n) \mathring\otimes (v_n+z_n) -\big(\chi_{\lambda_n}'(v)(u\mathring\otimes u)\big)*_x\phi_{\lambda_{n}^{-1}}*_t\psi_{\lambda_{n}^{-1}}-\lambda_{n}^{-2}\nabla(v_n+z_n)
	\\&\quad+\mathcal{R}\big((\chi_{\lambda_n}'(v)z)*_x\phi_{\lambda_{n}^{-1}}*_t\psi_{\lambda_{n}^{-1}}-z_n\big).
	\end{align*}
Here we take $\nu=\lambda_n^{-2}$.
To estimate $\$\mathring{R}_n\$_{L^1,r}$, we apply \eqref{b:z} and \eqref{bdvmom} to obtain
\begin{align*}
	&\$(v_n+z_n) \mathring\otimes (v_n+z_n) -\big(\chi_{\lambda_n}'(v)(u\,{\mathring\otimes} \, u)\big)*_x\phi_{\lambda_{n}^{-1}}*_t\psi_{\lambda_{n}^{-1}}\$_{L^1,r}
	\\&\lesssim\$v_n-v^{({\lambda_n})}\$_{L^2,2r}+\$v-v^{({\lambda_n})}\$_{L^2,2r}+\$z_n-z\$_{L^2,2r}
	\\&\quad+\$ \chi_{\lambda_n}'(v)-1\$_{L^q,\bar r}\$u\otimes u\$_{L^{p/2},r'}+\lambda_n^{-\bar\beta'}\$u\$_{H^{\bar\beta'},2r'}^2+\lambda_n^{-\bar\beta'}\$u\$_{C_t^{\bar\beta'}L^2,2r'}^2,
\end{align*}
\begin{align*}
	\$\lambda_{n}^{-2}\nabla(v_n+z_n)\$_{L^1,r}\lesssim \lambda_n^{-1}\$v_n\$_{L^2,2r}+\lambda_n^{-2}\lambda_{n+1}^{\alpha/8}\$z_n\$_{L^2,2r},
\end{align*}
and
\begin{align*}
	&\$\mathcal{R}\big((\chi_{\lambda_n}'(v)z)*_x\phi_{\lambda_{n}^{-1}}*_t\psi_{\lambda_{n}^{-1}}-z_n\big)\$_{L^1,r}
\\&\lesssim \$ \chi_{\lambda_n}'(v)-1\$_{L^q,\bar r}\$z\$_{L^2,r'}+(\lambda_{n}^{-\sigma}+\lambda_{n+1}^{-\alpha\sigma/8})\$z\$_{H^\sigma,r}+\lambda_{n}^{-1/2+\delta}\$z\$_{C_t^{1/2-\delta}L^2,r},
\end{align*}
where $\frac1{\bar r}+\frac1{r'}=\frac1r$, $\frac2p+\frac1q=1$.

Hence,  we use \eqref{bd:vM1}, \eqref{bd:vM2} \eqref{bdvn} and \eqref{bd:s5z1} to have
\begin{align*}
	\$\mathring{R}_n\$_{L^1,r}
	&\lesssim \lambda_{n}^{-1}+\lambda_{n}^{-(\bar\beta\wedge\sigma\wedge \bar\beta')}+\lambda_{n+1}^{-\alpha\sigma/8}+\lambda_{n}^{-1/2+\delta}\leq \frac1{48}\delta_{n+2}\underline{e},
\end{align*}
where we used  $\frac12-\frac1p=\bar\beta'/3$ to have $\|u\|_{L^p}\lesssim \|u\|_{H^{\beta'}}$ and we assumed $\alpha b/8<1$ for the parameters which is consistent with the parameters in Section \ref{s:par}.
For $m\geq1$
\begin{align*}
	\$\mathring{R}_n\$_{L^1,m}\leq 8mL^2+4(\lambda_n+1)^2+\lambda_n^{-1}(\lambda_n+1+2mL)+4mL\leq \lambda_n^2(24mL^2)^6.
\end{align*}
Here we may adjust $L$ to be large enough such that \eqref{b:z} holds and  to absorb extra constant. It is easy to see the bounds for $v_n$ in Proposition \ref{p:iteration} hold by choose $M_0$ large enough.
Now, we run the convex integration from Proposition \ref{p:iteration} with $\nu=\lambda_n^{-2}$ and we get a sequence $\{v_q\}_{q\geq n}$ and limit $v_{\nu} = \lim_{q\to\infty} v_q$ in $C(\mR,H^{\vartheta})\cap C^\vartheta(\mR,H^{\vartheta}) $ which
solves the stochastic Navier--Stokes equations with $\nu=\lambda_n^{-2}$. Moreover, using \eqref{eqtheta1} and \eqref{eq:app} it holds
\[ \$u - u_\nu \$_{L^2,r}\leq\$ u-v_n-z_n\$_{L^2,r}+ \sum_{q =
	n}^{\infty} \$(u_{q + 1} - u_q)\$_{L^2,r}\leq \eps.\]
Here in the last step we used \eqref{eqtheta1} and $a$ large enough.
The rest of the proof follows exactly the same arguments  as in the proof of Theorem \ref{thm:6.1}.
\end{proof}

Combining the above with the proof of Theorem \ref{th:5.1} we obtain the following result.

\bt\label{th:new1}Assume that $\tr((-\Delta)^\sigma GG^*)<\infty$ for some $\sigma>0$.
Let  $r>1$, $\varepsilon>0$, let $((\Omega,\cF,\bP),u, B)$ be a stationary solution to the stochastic Euler equations \eqref{2}  satisfying \eqref{eq:new} for some $\bar\beta'>0$, $r'>r$ and $L \geq (2\pi)^3$. There exist $\nu>0$, $\vartheta>0$ and stationary solutions $((\tilde\Omega,\tilde \cF, \tilde\bP),\tilde u,\tilde u_\nu,\tilde B)$ where $\mathcal{L}[\tilde u]=\mathcal{L}[u]$, $\tilde u$ solves \eqref{2} and $\tilde u_\nu$ solves  \eqref{1}  on $\mR\times \mT^3$, they belong to $ {C(\mR;H^{\vartheta})\cap C^{\vartheta}(\mR;L^{2})}$ $\mathbf{P}$-a.s.  and
\begin{align}\label{eq:00n}
\$\tilde u-\tilde u_\nu\$_{L^2,r}\leq \eps.
\end{align}
Furthermore, there exists a sequence of stationary solutions to the stochastic Navier--Stokes equations \eqref{1} with vanishing viscosities $\nu=\nu_n$ which converges in law in $C(\mathbb{R},L^2_\sigma)$ to $u$.
\et
\begin{proof}
The proof proceeds similarly to the proof of Theorem \ref{th:s1}. The main difference is that  now we consider the  ergodic averages
$$\nu_T:=\frac1{T}\int_0^T \cL[S_t(u, u_\nu, B)]\dif t,\qquad T\geq0.$$
The uniform bounds of $u,u_\nu$  imply tightness of $\nu_T$ in $C(\mR;L_{\sigma}^2)^3$.
Consequently, we obtain a stationary solution $((\tilde{\Omega},\tilde{\cF},\tilde{\bP}),\tilde{u},\tilde{u}_\nu,\tilde B)$ such that \eqref{eq:00n} holds.
Thus the result follows.
\end{proof}

\section{Stationary solutions to the deterministic Navier--Stokes/Euler equations}
\label{s:6}

In this section, we construct  random statistically stationary solutions to the deterministic Navier--Stokes/Euler equations on $\R\times\mathbb{T}^{3}$. As mentioned in the introduction, a lot  has been already achieved by using the known deterministic results about Euler and Navier--Stokes equations. Furthermore, also the results of Section~\ref{s:in}, Section~\ref{s:4} and Section~\ref{s:5} can be applied with $G=0$. Therefore, here we focus on proving that the constructed stationary solutions may be genuinely random as well as time dependent.  Precisely, we show  that  the solutions can be close in a certain sense to a given stationary stochastic process. This further highlights the fact how arbitrary the stationary solutions to the deterministic Navier--Stokes/Euler equations can be. In particular, the constructed stationary solutions can possess ``almost'' Gaussian or non-Gaussian statistics.

In the sequel, we make a few  modifications in the construction of Section~\ref{s:in}.
We consider the iteration
\begin{align}\label{iter:1}
\p_t u_q+\div (u_q\otimes u_q)+\nabla p_q=\nu\Delta u_q+\div \mathring{R}_q
\end{align}
with $\nu=1$ or $0$, which corresponds to the Navier--Stokes and Euler equations, respectively.
In this case Proposition \ref{p:iteration} holds for  $(u_q,\mathring{R}_{q})$ satisfying \eqref{iter:1}.
Its proof simplifies as we do not need to include the process $z_{q}$ anymore. Based on this, we obtain the following result.

\bt\label{thm:6.11}
Let $r>1$ be fixed and $Z$ be an $\mathcal{F}$-measurable stationary stochastic process  with smooth trajectories and vanishing
mean and divergence and satisfying
\begin{equation}\label{eq:Z}
\left(\mathbf{E}\|Z\|^m_{L^2}+\E\|Z\|^m_{C^2_{t,x}}\right)^{1/m}\leq m^{1/2}L,
\end{equation}
for any $m > 1$ and some $L \geq (2\pi)^3$. Let a smooth function $e:\mR\to(0,\infty)$ satisfying $\bar e\geq e(t)\geq \underline{e}>192 r L^2$ and $\varepsilon>0$ be given. There exists an $\mathcal{F}$-measurable  process $u$ which belongs to $ {C(\mR;H^{\vartheta})\cap C^{\vartheta}(\mR;L^{2})}$ $\mathbf{P}$-a.s.  and is an analytically weak solution to the deterministic Navier--Stokes/Euler equations   on $\mR\times \mT^3$.
Moreover, there exists $\vartheta>0$  such that
$$
\$u\$_{H^\vartheta,2r}+\$u\$_{C_{t}^\vartheta L^2,2r}<\infty,
$$
and for $t\in\mR$
\begin{align}
\label{eq:K1 1}
\begin{aligned}
	\mathbf{E}\|u(t)\|_{L^2}^2=e(t),%
\end{aligned}
\end{align}
and
\begin{align}\label{eq:uwp 1}
\sup_{t\in\R}\mathbf{E}\|u - Z\|_{C_{t}W^{1,1}}^r\leq \varepsilon.
\end{align}
\et

\begin{proof}
We start  the iteration with
$ u_0 = Z ,$
so \eqref{iter:1} on the level $q=0$ reads as
\[ \partial_t Z + \tmop{div} (Z \otimes Z) + \nabla p_0 - \nu\Delta Z=: \tmop{div} \mathring{R}_0 \]
with
\[ \mathring{R}_0 = Z \mathring\otimes Z -\mathcal{R}(\nu\Delta Z -\partial_t Z). \]
The last term is the reason why we require smoothness of trajectories of the process $Z$ and we also need a bound for the $C^2_{t,x}$-norm of $v_0=Z$.
Hence, for some $r>1$ such that $\underline{e}\geq 192 rL^2$
\begin{align*}
	\$\mathring{R}_0\$_{L^1,r}\leq 2rL^2+(2\pi)^3\cdot2rL\leq \frac1{48}\underline{e},
\end{align*}
and for $m\geq1$
\begin{align*}
	\$\mathring{R}_0\$_{L^1,m}\leq 2mL^2+(2\pi)^3\cdot2mL\leq 4mL^2.
\end{align*}
Now, we run the convex integration from Proposition \ref{p:iteration} and we get a limit $u = \lim_{q\to\infty} u_q$ in $C(\mR,H^{\vartheta})\cap C^\vartheta(\mR,H^{\vartheta}) $ which
solves the deterministic Navier--Stokes/Euler equations. Moreover, it holds
\[ u - Z =  \sum_{q =
	0}^{\infty} (u_{q + 1} - u_q).\]
Thus, 	as in Proposition \ref{p:iteration} we could choose $a$ large enough so that \eqref{bdvq:w1p} implies \eqref{eq:uwp 1}.
The rest of the proof follows exactly the same arguments  as in the proof of Theorem \ref{thm:6.1}.
\end{proof}

\begin{remark}\label{r:6.2}
From the proof it can be seen that the stochastic convex integration is not necessary provided $Z$ satisfies certain stronger assumptions. For instance, if $Z$ possesses a uniform in $\omega$ bound in {$C^{2}_{b}(\R\times\mT^{3})$},
a deterministic convex integration \`a la \cite{BV19a} can be applied pathwise, an $\omega$-dependent energy can be prescribed pathwise and the expectation in \eqref{eq:uwp 1} can be dropped.
Furthermore, the stochastic convex integration is also  not necessary provided the trajectories of $Z$ belong to {$C^{2}_{b}(\R\times\mT^{3})$} a.s. In this case, we can apply the deterministic convex integration on each of the sets
$$
\Omega_{L}=\left\{\omega\in\Omega;\, L-1\leq \|Z(\omega)\|_{{C^{2}_{b}(\R\times\mT^{3})}}< L\right\}
$$
and glue the solutions together similarly to \cite[Theorem 3.2]{HZZ22}. Restricting to the sets $\Omega_{L}$ permits to obtain the $\mathcal{F}$-measurability of the solutions, since  the parameters in the  convex integration only differ for different $L$ but not for each different  $\omega$.
\end{remark}

Combining the above with the proof of Theorem \ref{th:5.1} we obtain the following result.

\bt\label{th:6.21}
Let  $r>1$, $\varepsilon>0$, let $Z$ be as in Theorem \ref{thm:6.11} and let $K\geq192 rL^2$. There exist a random, time dependent, stationary solution $(\tilde\Omega,\tilde \cF,\tilde\bP,\tilde u)$ to the deterministic Navier--Stokes/Euler equations  on $\mR\times \mT^3$ satisfying
\begin{align*}
\begin{aligned}
	\tilde{\mathbf E}\|\tilde u\|_{L^2}^2=K,
\end{aligned}
\end{align*}
and a stochastic process $\tilde Z$   defined on the same probability space,   $\mathcal{L}[\tilde Z]=\mathcal{L}[Z]$, so that
\begin{align}\label{eq:6.44}
\tilde\E\|\tilde u-\tilde Z\|_{C_tW^{1,1}}^r\leq \eps.
\end{align} {Furthermore, there exist non-unique ergodic stationary solutions.}
\et
\begin{proof}
The proof proceeds similarly to the proof of Theorem \ref{th:s1}. The main difference is that  now we consider the  ergodic averages
$$\nu_T:=\frac1{T}\int_0^T \cL[S_t(u, Z)]\dif t:=\frac1{T}\int_0^T \cL[(u, Z)(t+\cdot)]\dif t,\qquad T\geq0.$$
The uniform bounds of $u$ and $Z$ imply tightness of $\nu_T$ in $C(\mR;L_{\sigma}^2)\times C(\mR;C^{1})$.
We  could find a stationary solution $(\tilde{\Omega},\tilde{F},\tilde{P},\tilde{u},\tilde{Z})$ such that $\tilde{\E}\|\tilde{u}\|_{L^2}^2=K$ for some given $K$ and
\begin{align*}
	\sup_{t\in\R}	\tilde{\E}\|\tilde{u}-\tilde{Z}\|_{C_tW^{1,{1}}}\leq \eps.
\end{align*}
This also implies that the stationary solution depends on time and is genuinely random provided $Z$ is time dependent and genuinely random. For different $K$ we have different stationary solutions which also gives non-unique ergodic stationary solutions.
\end{proof}

The above results offer a lot of freedom in the choice of the process $Z$, showing  how arbitrary the law of constructed stationary solutions can be. For instance, the construction can be  performed with a Gaussian stationary process $Z$ with smooth trajectories. A simple example is given by
$$
Z(t)=\left(\cos(t)\xi_{1}+\sin(t)\xi_{2}\right)e^{ik\cdot x}e_{k},
$$
where $\xi_{1}$, $\xi_{2}$ are independent standard Gaussians and $k, e_{k}\in\mR^{3}\setminus\{0\}$ are  orthogonal.
In this case, we could find $L>0$ such that for every $m\geq2$, $\$Z\$_{C^2_{t,x},m}+\$Z\$_{L^2,m}\leq (m-1)^{1/2}L$. Hence, \eqref{eq:Z} is satisfied.
Note that this corresponds to the second situation discussed in Remark~\ref{r:6.2}, namely, where the trajectories of $Z$ belong to {$C^{2}_{b}(\R\times\mT^{3})$} a.s.

To construct a Gaussian process $Z$ outside of the simplified setting of Remark~\ref{r:6.2}, we may consider $z$ as in Section~\ref{s:in}, i.e. a stationary solution to \eqref{linear}, and define $Z$ to be its space-time mollification. Then \eqref{eq:Z} follows from Proposition~\ref{fe z} and the stochastic convex integration as used in Theorem~\ref{thm:6.11} can be applied.

But we may as well  define $Z$ to be non-Gaussian. For example, we let
$
Z(t)=X(t)=\cos(t+Y)e^{ik\cdot x}e_k,
$
where $k,\, e_{k}\in\mathbb{R}^{3}\setminus\{0\}$ are orthogonal and $Y$ is uniformly
distributed on $(0, 2 \pi]$. Alternatively, we use the process $X$ given above and define $Z(t)=\int_0^\infty e^{s(\nu\Delta-I)}X(t-s)\dif s$ which is a stationary solution to the following equation
\begin{align}\label{eq:ZZ1}
\partial_t Z  + Z  = \nu\Delta Z+X,
\end{align}
which is also of zero mean and divergence free. Both these examples have uniform in $\omega$ bounds as discussed in Remark~\ref{r:6.2} and they are non-Gaussian.

\begin{remark}
Our solutions satisfy a weaker version of the so-called Kolmogorov hypothesis in the spirit of \cite{CG12,CV18} which is sufficient to perform rigorously the vanishing viscosity limit and obtain stationary solutions to the Euler equations.
\end{remark}

We conclude this section by an observation that stationary solutions to the deterministic Navier--Stokes/Euler equations can also be obtained as limits of stationary solutions to the stochastic counterparts of the equations, possibly combining with a vanishing viscosity limit.
More precisely, we have the following result.

\bt
\label{thm:6.4}
Suppose that $\tr((-\Delta)^{3/2+\sigma}GG^*)<\infty$ for some $\sigma>0$. Let  $r>1$ fixed, $Z$ be as in Theorem \ref{thm:6.11}, $K\geq384 rL^2$ and $\varepsilon>0$.
For an arbitrary sequence of vanishing constants $\gamma_{1,n}, \,\gamma_{2,n}\geq 0$,  $n\in\N$, $\gamma_{1,n}, \,\gamma_{2,n}\to 0$, up to a change of probability space, there exist a sequence of stationary solutions $u_n$, $n\in\N$, to the following stochastic Navier--Stokes/Euler equations (for $\nu=1$ or $\nu=0$)
\begin{align}\label{eq:unn:app}
\dif u_n+\div(u_n\otimes u_n)\,\dif t+\nabla P_{n}\,\dif t={(\nu+\gamma_{1,n})}\Delta u_n\,\dif t+{\gamma_{2,n}}\dif B,
\end{align}
so that the corresponding family of laws $\mathcal{L}[u_{n}]$, $n\in\N$, is tight in $C(\mR;L_{\sigma}^2)$ and
every accumulation point is a stationary solution to  the deterministic Navier--Stokes/Euler equations on $\mR\times \mT^3$. Furthermore, every accumulation point $u$ satisfies
\[\E\|u\|_{L^2}^2=K,\]
and
\begin{align}\label{u-Z}
\E\|u-Z\|_{C_tW^{1,1}}^r\leq \eps.
\end{align}\
\et

\begin{remark}
The formulation of the above theorem in particular permits a simultaneous limit of vanishing viscosity and vanishing noise. We recall that the particular choice of $\gamma_{2,n}=\sqrt{\gamma_{1,n}}$ was treated in \cite{K04} (see also \cite{GHSV15}) in two spatial dimensions. It gave rise to the so-called Kuksin measures, a genuinely random statistically stationary solutions to the deterministic Euler equations on $\mathbb{T}^{2}$. It was argued on page 472 in \cite{K04} that in three dimensions, the correct scaling is $\gamma_{2,n}=1$ in order to be consistent with Kolmogorov's prediction of anomalous dissipation, as suggested  by the formal energy equality.
\end{remark}

\begin{proof}[Proof of Theorem \ref{thm:6.4}]
We decompose \eqref{eq:unn:app} with $z_n$ a stationary solution to
\begin{align*}
	\dif z_n+z_n\dif t+\nabla P_{1,n}&=(\nu+\gamma_{1,n})\Delta z_n\dif t+\gamma_{2,n}\dif B,\\\div z_n&=0.
\end{align*}
Then we have for any $p\geq 2$, $0<\delta<1/2$
\begin{align*}
	\$z_n\$_{C_t^{1/2-\delta}L^\infty,p}+\$z_n\$_{C_t{H^{3/2}},p}\leq \gamma_{2,n}(p-1)^{1/2}L_1,
\end{align*}
for some $L_1>0$.
For each fixed $n$, we run the convex integration iteration as in Proposition \ref{p:iteration} indexed by $q$ and  starting from $v_0=Z$, which gives
\begin{align*}
	\$\mathring{R}_0\$_{L^1,r}\leq {4rL^2}+2\cdot(2\pi)^3rL+\gamma_{2,n}L_1+2\gamma^2_{2,n}rL^2_1\leq 8rL^2\leq\frac1{48}\underline{e}.
\end{align*}
Here we may choose $n$ large enough such that $\gamma_{2,n}4rL^2_1+\gamma_{2,n}L_1$ is small.
We then could use similar argument as in Theorem \ref{thm:6.11} and Theorem \ref{th:6.21} to construct, up to a change of probability space, stationary solutions $(u_n,z_n)$ so that
\begin{align*}
	\E\|u_n\|_{L^2}^2=K,
\end{align*}
and
\begin{align}\label{eq:unZ}
	\sup_{t\in\R}\E\|u_n-z_n-Z\|_{C_tW^{1,1}}^r\leq \eps.
\end{align}
Furthermore, for some $\vartheta>0$ we obtain
\begin{align*}
	\$u_n\$_{H^\vartheta,2r}+\$u_n\$_{C_t^\vartheta L^2,2r}\lesssim1,
\end{align*}
with the proportional constant independent of $n$.
Hence, we obtain tightness of $u_n$ in $C(\mR;L_{\sigma}^2)$ as in the proof of Theorem \ref{th:5.1} and conclude that the tight limit is a stationary solution to the deterministic Navier--Stokes/Euler equations. Since $z_n\to0$ in $C(\R;H^1)$, \eqref{eq:unZ} leads to \eqref{u-Z}.
\end{proof}

\appendix
\renewcommand{\appendixname}{Appendix~\Alph{section}}
\renewcommand{\theequation}{A.\arabic{equation}}

\section{Intermittent jets}
\label{s:B}

In this  part we recall the construction of  intermittent jets from \cite[Section 7.4]{BV19}.
We point out that the construction is entirely deterministic, that is, none of the functions below depends on $\omega$.
Let us begin with  the following geometric lemma which can be found in  \cite[Lemma 6.6]{BV19}.

\bl\label{geometric}
Denote by $\overline{B_{1/2}}(\mathrm{Id})$ the closed ball of radius $1/2$ around the identity matrix $\mathrm{Id}$, in the space of $3\times 3$ symmetric matrices. There
exists $\Lambda\subset \mathbb{S}^2\cap \mathbb{Q}^3$ such that for each $\xi\in \Lambda$ there exists a  $C^\infty$-function $\gamma_\xi:\overline{B_{1/2}}(\mathrm{Id})\rightarrow\mathbb{R}$ such that
\begin{equation*}
R=\sum_{\xi\in\Lambda}\gamma_\xi^2(R)(\xi\otimes \xi)
\end{equation*}
for every symmetric matrix satisfying $|R-\mathrm{Id}|\leq 1/2$.
For $C_\Lambda=8|\Lambda|(1+8\pi^3)^{1/2}$, where $|\Lambda|$ is the cardinality of the set $\Lambda$, we define
the constant
\begin{equation*}
M=C_\Lambda\sup_{\xi\in \Lambda}(\|\gamma_\xi\|_{C^0}+\sum_{|j|\leq N}\|D^j\gamma_\xi\|_{C^0}).
\end{equation*}
For each $\xi\in \Lambda$ let us define $A_\xi\in \mathbb{S}^2\cap \mathbb{Q}^3$ to be an orthogonal vector to $\xi$. Then for each $\xi\in\Lambda$ we have that $\{\xi, A_\xi, \xi\times A_\xi\}\subset \mathbb{S}^2\cap \mathbb{Q}^3$ form an orthonormal basis for $\mathbb{R}^3$.
We label by $n_*$ the smallest natural such that
\begin{equation*}\{n_*\xi, n_*A_\xi, n_*\xi\times A_\xi\}\subset \mathbb{Z}^3\end{equation*}
for every $\xi\in \Lambda$.
\el

Let $\Phi:\mathbb{R}^2\rightarrow\mathbb{R}$ be a smooth function with support in a ball of radius $1$. We normalize $\Phi$ such that
$\phi=-\Delta \Phi$ obeys
\begin{equation}\label{eq:phi}
\frac{1}{4\pi^2}\int_{\mathbb{R}^2}\phi^2(x_1,x_2)\dif x_1\dif x_2=1.
\end{equation}
By definition we know $\int_{\mathbb{R}^2}\phi dx=0$. Define $\psi:\mathbb{R}\rightarrow\mathbb{R}$ to be a smooth, mean zero function with support in the ball of radius $1$ satisfying
\begin{equation}\label{eq:psi}
\frac{1}{2\pi}\int_{\mathbb{R}}\psi^2(x_3)\dif x_3=1.
\end{equation}
For parameters $r_\perp, r_\|>0$ such that
\begin{equation*}r_\perp\ll r_\|\ll1,\end{equation*}
we define the rescaled cut-off functions
\begin{equation*}\phi_{r_\perp}(x_1,x_2)=\frac{1}{r_\perp}\phi\left(\frac{x_1}{r_\perp},\frac{x_2}{r_\perp}\right),\quad
\Phi_{r_\perp}(x_1,x_2)=\frac{1}{r_\perp}\Phi\left(\frac{x_1}{r_\perp},\frac{x_2}{r_\perp}\right),\quad \psi_{r_\|}(x_3)=\frac{1}{r_\|^{1/2}}\psi\left(\frac{x_3}{r_\|}\right).\end{equation*}
We periodize $\phi_{r_\perp}, \Phi_{r_\perp}$ and $\psi_{r_\|}$ so that they are viewed as periodic functions on $\mathbb{T}^2, \mathbb{T}^2$ and $\mathbb{T}$ respectively.

Consider a large real number $\lambda$ such that $\lambda r_\perp\in\mathbb{N}$, and a large time oscillation parameter $\mu>0$. For every $\xi\in \Lambda$ we introduce
\begin{equation*}\aligned
\psi_{(\xi)}(t,x)&:=\psi_{\xi,r_\perp,r_\|,\lambda,\mu}(t,x):=\psi_{r_{\|}}(n_*r_\perp\lambda(x\cdot \xi+\mu t))
\\ \Phi_{(\xi)}(x)&:=\Phi_{\xi,r_\perp,\lambda}(x):=\Phi_{r_{\perp}}(n_*r_\perp\lambda(x-\alpha_\xi)\cdot A_\xi, n_*r_\perp\lambda(x-\alpha_\xi)\cdot(\xi\times A_\xi))\\
\phi_{(\xi)}(x)&:=\phi_{\xi,r_\perp,\lambda}(x):=\phi_{r_{\perp}}(n_*r_\perp\lambda(x-\alpha_\xi)\cdot A_\xi, n_*r_\perp\lambda(x-\alpha_\xi)\cdot(\xi\times A_\xi)),
\endaligned\end{equation*}
where $\alpha_\xi\in\mathbb{R}^3$ are shifts to ensure that $\{\Phi_{(\xi)}\}_{\xi\in\Lambda}$ have mutually disjoint support.

The intermittent jets $W_{(\xi)}:\mathbb{R}\times\mathbb{T}^3 \rightarrow\mathbb{R}^3$ are defined as in \cite[Section 7.4]{BV19}.
\begin{equation}\label{intermittent}W_{(\xi)}(t,x):=W_{\xi,r_\perp,r_\|,\lambda,\mu}(t,x):=\xi\psi_{(\xi)}(t,x)\phi_{(\xi)}(x).\end{equation}
By the choice of $\alpha_\xi$ we have that
\begin{equation}\label{Wxi}
W_{(\xi)}\otimes W_{(\xi')}\equiv0, \textrm{ for } \xi\neq \xi'\in\Lambda,
\end{equation}
and by the normalizations \eqref{eq:phi} and \eqref{eq:psi} we obtain
$$
\frac1{(2\pi)^3}\int_{\mathbb{T}^3}W_{(\xi)}(t,x)\otimes W_{(\xi)}(t,x)\dif x=\xi\otimes\xi.
$$
These facts combined with Lemma \ref{geometric} imply that
\begin{equation}\label{geometric equality}
\frac1{(2\pi)^3}\sum_{\xi\in\Lambda}\gamma_\xi^2(R)\int_{\mathbb{T}^3}W_{(\xi)}(t,x)\otimes W_{(\xi)}(t,x)\dif x=R,
\end{equation}
for every symmetric matrix $R$ satisfying $|R-\textrm{Id}|\leq 1/2$. Since $W_{(\xi)}$ are not divergence free, we  introduce the corrector term
\begin{equation}\label{corrector}
W_{(\xi)}^{(c)}:=\frac{1}{n_*^2\lambda^2}\nabla \psi_{(\xi)}\times \textrm{curl}(\Phi_{(\xi)}\xi)
=\textrm{curl\,curl\,} V_{(\xi)}-W_{(\xi)}.
\end{equation}
with
\begin{equation*}
V_{(\xi)}(t,x):=\frac{1}{n_*^2\lambda^2}\xi\psi_{(\xi)}(t,x)\Phi_{(\xi)}(x).
\end{equation*}
Thus we have
\begin{equation*}
\div\left(W_{(\xi)}+W_{(\xi)}^{(c)}\right)\equiv0.
\end{equation*}

Finally, we  recall the key   bounds from \cite[Section 7.4]{BV19}. For $N, M\geq0$ and $p\in [1,\infty]$ the following holds
\begin{equation}\label{bounds}\aligned
&\|\nabla^N\partial_t^M\psi_{(\xi)}\|_{C_{t}L^p}\lesssim r^{1/p-1/2}_\|\left(\frac{r_\perp\lambda}{r_\|}\right)^N
\left(\frac{r_\perp\lambda \mu}{r_\|}\right)^M,\\
&\|\nabla^N\phi_{(\xi)}\|_{L^p}+\|\nabla^N\Phi_{(\xi)}\|_{L^p}\lesssim r^{2/p-1}_\perp\lambda^N,\\
&\|\nabla^N\partial_t^MW_{(\xi)}\|_{C_{t}L^p}+\frac{r_\|}{r_\perp}\|\nabla^N\partial_t^MW_{(\xi)}^{(c)}\|_{C_{t}L^p}+\lambda^2\|\nabla^N\partial_t^MV_{(\xi)}
\|_{C_{t}L^p}\\&\lesssim r^{2/p-1}_\perp r^{1/p-1/2}_\|\lambda^N\left(\frac{r_\perp\lambda\mu}{r_\|}\right)^M,
\endaligned\end{equation}
where the implicit constants may depend on $p, N$ and $M$, but are independent of $\lambda, r_\perp, r_\|, \mu$.

\renewcommand{\theequation}{B.\arabic{equation}}

\section{Estimates of $\rho$ and $a_{(\xi)}$}\label{ap:B}

For completeness, we include here the detailed proof  of the estimates \eqref{estimate aN} and \eqref{estimate aN0} employed in Section~\ref{s:313}.

We first aim at estimating the $C^{N}_{t,x}$-norm of $\rho$ for $N\in\mathbb{N}$. To this end, we first apply  the chain rule  \cite[Proposition C.1]{BDLIS16} to the function $\Psi(z)=\sqrt{\ell^{2}+z^{2}}$, $|D^{m}\Psi(z)|\lesssim \ell^{-m+1}$ to  obtain
\begin{align}\label{rhoN1}
\|\rho\|_{C^N_{t,x}} \lesssim \left\| \sqrt{\ell^2 + | \mathring{R}_\ell |^2} \right\|_{C^{0}_{t,x}}+\|D \Psi \|_{C^{0}}
\| \mathring{R}_\ell \|_{C^N_{t,x}} + \| D\Psi \|_{C^{N-1}} \|
\mathring{R}_\ell \|^N_{C^1_{t,x}}+\|\gamma_{\ell}\|_{C^N_t} \\
\lesssim \ell^{- 4 - N} \|\mathring{R}_q\|_{C_{[t-1,t+1]}L^1} + \ell^{- 6N + 1}
\|\mathring{R}_q\|_{C_{[t-1,t+1]}L^1}^N+\ell^{-N}\delta_{q+1}\bar e. \nonumber\end{align}

For the amplitude functions $a_{(\xi)}$ defined in \eqref{amplitudes} we deduce using (\ref{rho})
\begin{equation}\label{estimate a1}
\begin{aligned}
	\|a_{(\xi)}\|_{C_tL^2}&\leq \|\rho\|_{C_tL^1}^{1/2}\|\gamma_\xi\|_{C^0(B_{1/2}(\Id))}\leq \frac{M}{8|\Lambda|(1+8\pi^{3})^{1/2}}\left(2\ell(2\pi)^{3}+2\|\mathring{R}_q\|_{C_{[t-1,t+1]}L^1}+\frac12\delta_{q+1}\bar e\right)^{1/2}
	\\&\leq\frac{M}{4|\Lambda|}\left(2\|\mathring{R}_q\|_{C_{[t-1,t+1]}L^1}+\frac12\delta_{q+1}\bar e\right)^{1/2},
\end{aligned}
\end{equation}
where   $M$ denotes the universal constant from Lemma~\ref{geometric}.

Let us now estimate the $C^{N}_{t,x}$-norm of $a_{(\xi)}$. By Leibniz rule, we get
\begin{equation}\label{aa}
\|a_{(\xi)}\|_{C^{N}_{t,x}}\lesssim\sum_{m=0}^N\|\rho^{\frac12}\|_{C^m_{t,x}}\left\|\gamma_{\xi}\left(\Id
-\frac{\mathring{R}_\ell}{\rho}\right)\right\|_{C^{N-m}_{t,x}}
\end{equation}
and estimate each norm separately.
First, by \eqref{rho0}
\begin{equation*}
\| \rho^{1 / 2} \|_{C^0_{t,x}} \lesssim \ell^{-2}\|\mathring{R}_q\|_{C_{[t-1,t+1]}L^1}^{1/2}+\delta_{q+1}^{1/2}\bar e^{1/2},
\end{equation*}
and by Lemma~\ref{geometric}
$$
\left\|\gamma_{\xi}\left(\Id
-\frac{\mathring{R}_\ell}{\rho}\right)\right\|_{C^{0}_{t,x}}\lesssim 1.
$$
Second, applying \cite[Proposition C.1]{BDLIS16} to the function
$
\Psi (z) = z^{1 / 2},$ $  | D^{m}\Psi(z) | \lesssim |z|^{1 / 2 - m}, $ for $m=1,\dots, N$, and using \eqref{rhoN1} and $\rho\geq \ell$ and $\ell^{-1}\geq \bar e$ we obtain for $m\geq1$
\begin{equation}\label{eq:rho12}
\begin{aligned}
	\| \rho^{1 / 2} \|_{C^m_{t,x}} &\lesssim \| \rho^{1 / 2} \|_{C^0_{t,x}}+\ell^{- 1 / 2} \|\rho\|_{C^{m}_{t,x}} + \ell^{1 / 2 -
		m} \|\rho\|_{C^{1}_{t,x}}^m\\
	& \lesssim \ell^{-3m+1/2}\delta_{q+1}^m+\ell^{- 2} \|\mathring{R}_q\|_{C_{[t-1,t+1]}L^1}^{1/2} +
	\ell^{-9 / 2 - m}  \|\mathring{R}_q\|_{C_{[t-1,t+1]}L^1}\\&+\ell^{-6 m+1/2}  \|\mathring{R}_q\|^m_{C_{[t-1,t+1]}L^1}
	+\ell^{-m-{3/2}}\delta_{q+1}+\ell^{-1/2}\delta_{q+1}^{1/2}.
\end{aligned}
\end{equation}
For $m\geq1$ using $\delta_{q+1}\leq 1$ the above is bounded by
$$\ell^{-6 m+1/2}(1+  \|\mathring{R}_q\|_{C_{[t-1,t+1]}L^1}^m).$$

We proceed with a bound for $\left\|\gamma_{\xi}\left(\Id
-\frac{\mathring{R}_\ell}{\rho}\right)\right\|_{C^{N-m}_{t,x}}$ for $m=0,\dots,N-1$. Keeping \cite[Proposition C.1]{BDLIS16} as well as Lemma~\ref{geometric} in mind, we need to estimate
\begin{equation}\label{gamma}
\left\| \frac{\mathring{R}_\ell}{\rho} \right\|_{C^{N - m}_{t,x}} + \left\| \frac{\nabla_{t,x}\mathring{R}_\ell}{\rho}
\right\|_{C^0_{t,x}}^{N - m} + \left\| \frac{\mathring{R}_\ell}{\rho^{2}} \right\|_{C^0_{t,x}}^{N - m} \|
\rho \|_{C^1_{t,x}}^{N - m} .
\end{equation}
We use $\rho\geq \ell$ to have
$$ \left\| \frac{\nabla_{t,x}\mathring{R}_\ell}{\rho} \right\|_{C^0_{t,x}}^{N - m} \lesssim \ell^{- (N - m)}
\ell^{(- 4 - 1) (N - m)}  \|\mathring{R}_q\|^{N-m}_{C_{[t-1,t+1]}L^1}\lesssim \ell^{-
6 (N - m)}\|\mathring{R}_q\|^{N-m}_{C_{[t-1,t+1]}L^1} , $$
and	in view of $|\frac{\mathring{R}_\ell}{\rho}|\leq 1$
$$ \left\| \frac{\mathring{R}_\ell}{\rho^2} \right\|_{C^0_{t,x}}^{N - m} \lesssim \left\| \frac{1}{\rho} \right\|_{C^0_{t,x}}^{N - m}\lesssim \ell^{-(N-m)},$$
and by \eqref{rhoN1} and $\bar e\delta_{q+1}\leq \ell^{-1}$
$$ \| \rho \|_{C^1_{t,x}}^{N - m} \lesssim \ell^{-5 (N - m)} \|\mathring{R}_q\|^{N-m}_{C_{[t-1,t+1]}L^1}+{\ell^{-2(N-m)}}. $$

To estimate the  first term in \eqref{gamma}, we write		
\begin{equation}\label{RR}
\left\| \frac{\mathring{R}_\ell}{\rho} \right\|_{C^{N - m}_{t,x}} \lesssim \sum_{k = 0}^{N - m}
\| \mathring{R}_\ell \|_{C^k_{t,x}} \left\| \frac{1}{\rho} \right\|_{C^{N - m - k}_{t,x}},
\end{equation}
where for $N-m-k=0$ we have
$$
\left\| \frac{1}{\rho} \right\|_{C^{0}_{t,x}} \lesssim \ell^{-1}
$$
and for $k=0,\dots,N-m-1$ using \eqref{rhoN1}
$$
\left\| \frac{1}{\rho} \right\|_{C^{N - m - k}_{t,x}} \lesssim \left\|\frac1\rho\right\|_{C^{0}_{t,x}}+\ell^{- 2} \| \rho
\|_{C^{N - m - k}_{t,x}} + \ell^{- (N - m - k){-1}} \| \rho \|_{C^1_{t,x}}^{N - m - k} $$
$$
\lesssim \ell^{{ - 4(N - m - k)-1}} + \ell^{-6- (N - m - k)}\|\mathring{R}_q\|_{C_{[t-1,t+1]}L^1}+
\ell^{- 6 (N - m - k)-1}\|\mathring{R}_q\|^{N-m-k}_{C_{[t-1,t+1]}L^1} . $$
Altogether, we therefore bound \eqref{RR} as
\begin{equation}\label{RR2}
\begin{aligned}
	\left\| \frac{\mathring{R}_\ell}{\rho} \right\|_{C^{N - m}_{t,x}} &\lesssim
	\ell^{-5-4(N-m)}\|\mathring{R}_q\|_{C_{[t-1,t+1]}L^1}+\ell^{-10-(N-m)}\|\mathring{R}_q\|_{C_{[t-1,t+1]}L^1}^2\\&+\ell^{-5-6(N-m)}{(\|\mathring{R}_q\|_{C_{[t-1,t+1]}L^1}^{N-m\rmb{+1}}+1)}.
\end{aligned}
\end{equation}

Finally, plugging \eqref{RR2}, and the other bounds into  \eqref{gamma} leads to
\begin{align*}
\left\|\gamma_{\xi}\left(\Id
-\frac{\mathring{R}_\ell}{\rho}\right)\right\|_{C^{N-m}_{t,x}}\lesssim& \ell^{-5-4(N-m)}\|\mathring{R}_q\|_{C_{[t-1,t+1]}L^1}+\ell^{-10-(N-m)}\|\mathring{R}_q\|_{C_{[t-1,t+1]}L^1}^2\\&+\ell^{-5-6(N-m)}(\|\mathring{R}_q\|_{C_{[t-1,t+1]}L^1}^{N-m\rmb{+1}}{+1})+\ell^{-{3}(N-m)}.
\end{align*}
For $N-m\geq1$ the above is bounded by
$$\ell^{-5-6(N-m)}(\|\mathring{R}_q\|^{N-m\rmb{+1}}_{C_{[t-1,t+1]}L^1}+\|\mathring{R}_q\|^{2}_{C_{[t-1,t+1]}L^1}+1).$$
Combining this
with the bounds for $\rho^{1/2}$ above and plugging into \eqref{aa} yields for $N\geq 2$
\begin{equation*}
\begin{aligned}
	\|a_{(\xi)}\|_{C^N_{t,x}}& \lesssim \ell^{-7-6N}(\|\mathring{R}_q\|_{C_{[t-1,t+1]}L^1}+1)^{N+\rmb{3/2}}.
\end{aligned}
\end{equation*}
It is easy to see that
\begin{equation}\label{estimate aN01}
\begin{aligned}
	\|a_{(\xi)}\|_{C^0_{t,x}}& \lesssim \ell^{-2}(\|\mathring{R}_q\|_{C_{[t-1,t+1]}L^1}+1)^{1/2}.
\end{aligned}
\end{equation}
By interpolation \rmb{$\|a_{(\xi)}\|_{C^N_{t,x}}\lesssim \|a_{(\xi)}\|_{C^0_{t,x}}^{1/2} \|a_{(\xi)}\|_{C^{2N}_{t,x}}^{1/2}$}, the following estimate
\begin{equation}\label{estimate aN1}
	\begin{aligned}
		\|a_{(\xi)}\|_{C^N_{t,x}}& \lesssim \ell^{-7-6N}(\|\mathring{R}_q\|_{C_{[t-1,t+1]}L^1}+1)^{N+\rmb{1}}.
	\end{aligned}
\end{equation}
 holds.
 
 \section*{Declarations}
 
 The authors have no competing interests to declare that are relevant to the content of this article.
 
 Data sharing not applicable to this article as no datasets were generated or analyzed during the current study.

%

\def\cprime{$'$} \def\ocirc#1{\ifmmode\setbox0=\hbox{$#1$}\dimen0=\ht0
  \advance\dimen0 by1pt\rlap{\hbox to\wd0{\hss\raise\dimen0
  \hbox{\hskip.2em$\scriptscriptstyle\circ$}\hss}}#1\else {\accent"17 #1}\fi}

\end{document}